\documentclass[11pt]{article}
\usepackage[dvips]{graphicx}
\usepackage{amssymb}

\def\R{\mathbf{R}}
\def\Z{\mathbf{Z}}
\def\e{\epsilon}

\def\con{\mathrm{const.}}

\def\proof{\par\medskip\noindent {\bf Proof.}{\hskip1em}}
\def\qed{~\vrule height .9ex width .8ex depth -.1ex}
\def\n#1{\|\hskip-1pt #1\hskip-1pt\|}
\def\f{\bar{f}}
\def\vol{\mathrm{vol}}
\def\lng{\mathrm{length}}
\def\area{\mathrm{area}}

\newtheorem{thm}{Theorem}
\newtheorem{lem}{Lemma}
\newtheorem{cor}[lem]{Corollary}
\newtheorem{prop}[lem]{Proposition}
\newtheorem{defi}[lem]{Definition}
\newtheorem{rem}[lem]{Remark}
\newtheorem{ques}{Question}

\title{Calibrations and isoperimetric profiles
\footnote{
\eject\noindent Work done within the GNSAGA group of INDAM, supported by MURST of Italy and by the Galileo project \emph{Asymptotic properties of groups and manifolds}. }}
\author{Renata Grimaldi\footnote{Univ Palermo,
Dipartimento di Metodi e Modelli Matematici, Palermo, I-90128.} and Pierre Pansu$^{1,2}$\footnote{$^{1}$ Univ Paris-Sud, Laboratoire de Math\'ematiques d'Orsay, Orsay, F-91405 ;\hfill\eject\indent\hskip7.8pt $^{2}$ CNRS, Orsay, F-91405.}}
\begin{document}

\maketitle

\begin{quote}
{\small 
ABSTRACT: We equip many non compact non simply connected surfaces with smooth Riemannian metrics whose isoperimetric profile is smooth, a highly non generic property. The computation of the profile is based on a calibration argument, a rearrangement argument, the Bol-Fiala curvature dependent inequality, together with new results on the profile of surfaces of revolution and some hardware know-how.
\par\medskip
RESUM\'E: On construit, sur des surfaces non compactes non simplement connexes, des m\'etriques riemanniennes lisses dont le profil isop\'erim\'etrique est lisse, ph\'enom\`ene hautement non g\'en\'erique. Le calcul du profil est bas\'e sur l'id\'ee de calibration, un argument de r\'earrangement, l'in\'egalit\'e de Bol-Fiala, un r\'esultat nouveau sur le profil des surfaces de r\'evolution, et des connaissances en plomberie.
}
\end{quote}

\section{Introduction}

\subsection{The problem}

Let $M$ be a smooth manifold. We are concerned with \emph{isoperimetric profiles} of metrics on $M$.

Fix a Riemannian metric on $M$ with total volume $V$. Given $v<V$, consider all domains (relatively compact open subsets with smooth boundary) in $M$ with volume $v$. Define $I_{M}(v)$ as the least upper bound of the boundary volumes of such domains. In this way, one gets a function $I_{M}:(0,V)\to\R_{+}$ called the \emph{isoperimetric profile} of $M$. 

For instance, Euclidean space $\mathbf{ R}^{n}$ has
\begin{displaymath}
I_{\mathbf{R}^{n}}(v)=\con\,v^{\frac{n-1}{n}}
\end{displaymath}
for some constant $\con(n)$. Indeed, for each $v$, round balls of volume $v$ minimize boundary volume among all domains of volume $v$. 

Here is an other example, where the isoperimetric profile is much easier to compute. Let $\Gamma$ be a lattice of parabolic translations of  hyperbolic $n$-space $\mathbf{ H}^{n}$ (i.e. $\Gamma$ has a unique fixed point $\zeta$ on the ideal boundary, and acts as a lattice of translations on the complement of $\zeta$, identified with euclidean $n-1$-space). Let $M=\mathbf{ H}^{n}/\Gamma$. We call such a manifold a \emph{complete constant curvature cusp}. Then
\begin{displaymath}
I_{M}(v)=(n-1)v.
\end{displaymath}
This will be proven as Corollary \ref{cusp} in section \ref{cali}.

In both examples, the isoperimetric profile $I$ is a smooth function. This is rather exceptional and related with the topological simplicity of the underlying manifolds. Smoothness of the isoperimetric profile is related to uniqueness of minimizers in the following way. For a generic Riemannian manifold, one expects that domains whose boundary has constant mean curvature come in finitely many smooth families. To each such family $D_{i,t}$, there corresponds a function $P_{i}$ defined by
\begin{displaymath}
P_{i}(vol(D_{i,t}))=vol(\partial D_{i,t})
\end{displaymath}
on some interval. The isoperimetric profile of $M$ is the minimum of these functions $P_{i}$. If $M$ has nontrivial topology, there must be at least two different families, and it seems very unlikely that the isoperimetric profile be everywhere smooth. This suggests the following question.

\begin{ques}\label{question}
Does every manifold admit a smooth (resp. real analytic) metric with smooth (resp. real analytic) isoperimetric profile ?
\end{ques}

\subsection{The results}

In this note, we construct examples of smooth complete metrics with smooth isoperimetric profiles on topologically non trivial manifolds, mostly in dimension 2. 

\begin{thm}
\label{nonachieved}
Every smooth (resp. real analytic) non compact manifold admits smooth (resp. real analytic) complete metrics of infinite volume whose isoperimetric profile vanishes identically.
\end{thm}

For such metrics, there are no extremal domains. That's cheating! Say a Riemannian manifold $M$ has an \emph{achieved} isoperimetric profile if for each $v\in(0,\vol(M))$, there exists a domain $D$ with $\vol(D)=v$ and $\vol(\partial D)=I_M (v)$. The isoperimetric profile of a compact Riemannian manifold is always achieved. Here is a simple sufficient condition for non compact surfaces.

\begin{defi}
\label{negative}
Say a non compact 2-dimensional Riemannian manifold $M$ is \emph{ultrahyperbolic} if
\begin{itemize}
  \item $M$ is complete;
  \item its curvature tends to $-\infty$ at infinity;
  \item its injectivity radius tends to $+\infty$ at infinity.
\end{itemize}
\end{defi}
On such a surface, the isoperimetric profile is achieved, this is proven as Proposition \ref{achieved} in the appendix.

\begin{thm}
\label{main}
Let $M$ be an orientable 2-dimensional manifold. Assume that one of the following properties holds.
\begin{itemize}
  \item $M$ has at least 4 ends.
  \item $M$ has 3 ends, at least one of which having infinite genus.
  \item $M$ has 2 ends and both have infinite genus. 
\end{itemize}
Then $M$ admits a complete smooth Riemannian metric of infinite area with \emph{achieved} isoperimetric profile $I_{M}(a)=a$ for all $a\in(0,+\infty)$.
\end{thm}

An even larger family of surfaces admit smooth metrics with smooth isoperimetric profiles.

\begin{thm}
\label{2ends}
Let $M$ be an orientable 2-dimensional manifold. Assume that one of the following properties holds.
\begin{itemize}
  \item $M$ has at least 2 ends.
  \item $M$ has one end of infinite genus.
\end{itemize}
Then $M$ admits a complete smooth Riemannian metric of infinite area with achieved isoperimetric profile whose square $I_{M}^2$ extends to a smooth function on some neighborhood of $\R_+$. If furthermore $M$ has finite topological type, one can arrange that $M$ be ultrahyperbolic.
\end{thm}

\begin{ques}
\label{quesconf}
Given a manifold $M$, can one construct metrics with smooth isoperimetric profile in all conformal classes on $M$ ? 
\end{ques}

\subsection{Scheme of proof}

The proof of Theorem \ref{nonachieved} is straightforward : given a function with compact sublevel sets on $M$, one arranges so that the volume of level sets tends to 0 although the total volume is infinite.

The metrics of Theorem \ref{main} are modelled on the 2-dimensional complete constant curvature $-1$ cusp $M_{0}$, where the solutions of the isoperimetric problem are the sublevels $V=\{f\leq t\}$ of a function $f$ such that $\Delta f=1$. In $M_{0}$, each level set $\partial V=\{f=t\}$ is \emph{calibrated} by the unit 1-form $\displaystyle \omega=*df$, i.e. 
\begin{displaymath}
\omega_{|\partial V}=vol_{\partial V}.
\end{displaymath}

The calibration method of Harvey and Lawson \cite{HarveyLawson} can be readily adapted to the isoperimetric problem, see section \ref{cali}. Given a calibrating 1-form (not quite a calibration since it is not closed), the method yields the value of the isoperimetric profile $I(v)$ at $v$ provided there exists a set with area $v$ and calibrated boundary.

For the surfaces with complicated topology dealt with in Theorem \ref{main}, the construction goes in two steps. 

First, glue together pieces of complete constant curvature $-1$ cusps in order to obtain a \emph{singular} surface with constant curvature $-1$, equipped with a \emph{level function} $f$ and a calibrating 1-form $*df$. This is done in section \ref{singular}.
This kind of surface has plenty of subsets with calibrated boundary. 
\begin{figure}
\begin{center}
\includegraphics[width=2in]{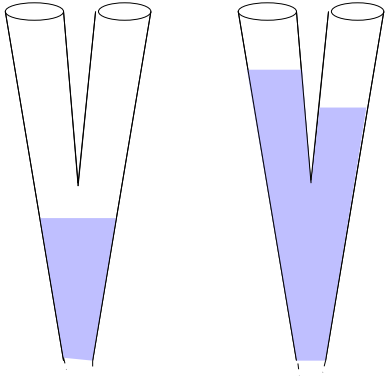}
\caption{Calibrated regions in a singular surface}
\label{calibres}
\end{center}
\end{figure}
Figure \ref{calibres} shows an example of a triply punctured sphere made of 3 pieces of a cusp. It has 2 families of calibrated subsets exhibiting different topological types. One family depends on one parameter, the second on two parameters. 

\begin{figure}
\begin{center}
\includegraphics[width=2in]{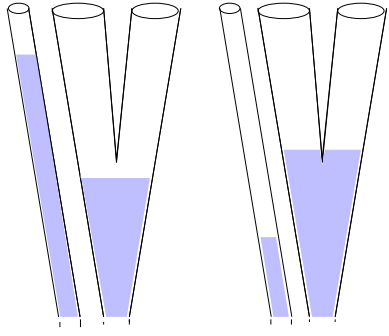}
\caption{The pipe clearing trick}
\label{calibrebis}
\end{center}
\end{figure}

Second, smooth away singularities, while maintaining a calibrating 1-form. This is done in section \ref{smoothing}. The resulting surface has less calibrated subsets, since levels passing close to singular points are not calibrated any more. Sometimes, they are still sufficient for the determination of the isoperimetric profile. This not the case for the example of figure \ref{calibres}. Indeed, there is a gap, around the area of the critical sublevel, in the values of the area which are achieved by calibrated subsets. 

To overcome this difficulty, the ''pipe clearing trick'' is applied. Imagine the surface as a pipestry and the calibrated regions as partial fillings of the pipestry with water. View a singular point as an obstacle blocking a pipe A. The traditional plumber's strategy consists in letting the water level raise in a neighboring pipe B, then pushing energeticly the water out of B so that in A water traverses the critical level, see figure \ref{calibrebis}. This allows, without change in the total amount of water, a change in the topology of the wet region. This hardware knowledge is turned into mathematics in section \ref{proof}. This trick is related to Pitts' observation that minimax cycles on a surface cannot have too many selfintersections, compare the four-legged star fish in \cite{Pitts} and \cite{Calabi-Cao} page 540.

The existence of a neighboring pipe is related to a Morse function (water level) having disconnected level sets. Such functions exist only on surfaces with many ends, as shown in section \ref{morse}.

\medskip

The examples of Theorem \ref{2ends} require two more steps.

When applied to a manifold $M$ with two points removed, Theorem \ref{main} yields a manifold with two constant curvature $-1$ cusps. Closing the cusps with suitably chosen caps of revolution yields a metric on $M$ whose isoperimetric profile is given by that of the caps, see section \ref{filling}. The only surfaces of revolution whose isoperimetric profiles were previously known to be smooth, see \cite{BC}, \cite{HHM}, \cite{R}, \cite{T} are those whose curvature is a nonincreasing function of the distance to the pole. For our purposes, this class needs to be slightly enlarged. An argument based on the strict stability of parallels shows that the property that rotationally symmetric domains are extremal is stable under small perturbations, see section \ref{revolution}. This provides us with the suitable caps.
\begin{figure}
\begin{center}
\includegraphics[width=4in]{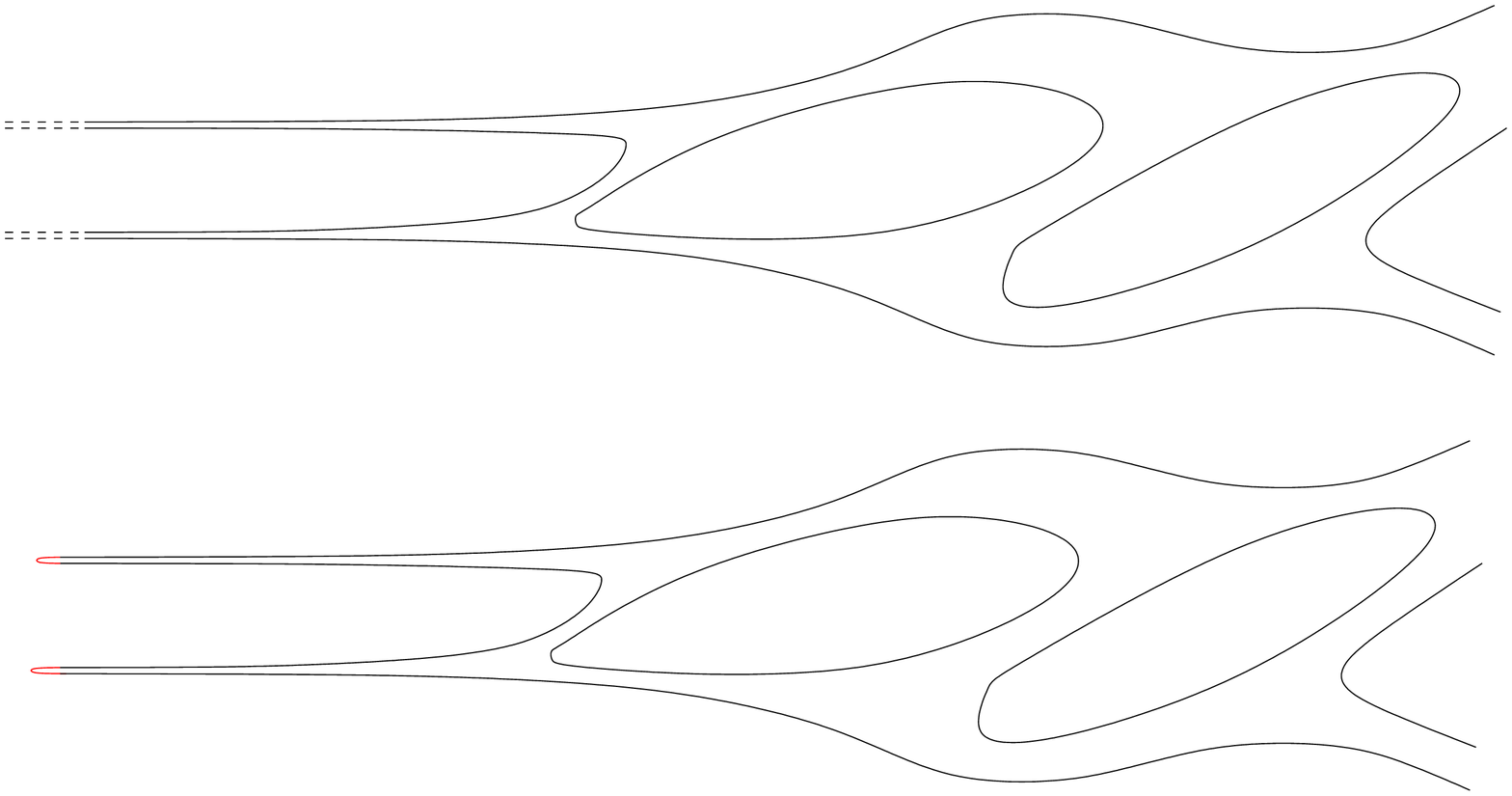}
\caption{Cusp closing}
\label{caps}
\end{center}
\end{figure}

Finally, a rearrangement argument shows that the calibration method extends to forms $\omega$ such that $d\omega=u\,vol$ where extremal domains are sublevel sets of $u$. Therefore an extra conformal change produces metrics whose isoperimetric profile equals any prescribed smooth convex function. Such metrics can be arranged to have curvature tending to $-\infty$ and injectivity radius tending to $+\infty$, see section \ref{rearrangement}.

\subsection{Acknowledgements}
This work owes a lot to V. A. Zorich and V. M. Kesel'man's construction of a singular metric with linear isoperimetric profile in a conformal class, \cite{KZ1}. It  arose from an attempt to smooth their metric. Many thanks are due to the anonymous referee for his help in improving the exposition.

\section{Non achieved isoperimetric profiles}

We prove Theorem \ref{nonachieved}.

Let $M$ be a smooth manifold, $w:M\to\R_+$ a smooth Morse function with compact sublevel sets and isolated critical values. Start with any complete smooth metric $g$ on $M$. Consider the positive continuous function $v$ defined on $(0,+\infty)$ by
\begin{eqnarray*}
v(t)=\vol(\{w\leq t\}).
\end{eqnarray*}
Let $u:\R+\to(0,+\infty)$ be a smooth function such that $uv$ tends to 0. In the conformal metric $g'=u(w)^{2/(n-1)}g$, the level set $\{w=t\}$ has volume $uv(t)$. Let $h:\R+\to\R_+$ be a smooth function such that $\displaystyle\int_{1}^{+\infty}huv(t)\,dt=+\infty$. Let $g''=g'+h(w)^2\,dw^2$. Then the volume of level sets of $w$ still tends to 0 since it does not change. Furthermore, for all $0<s<t$,
\begin{eqnarray*}
\vol''(\{s<w<t\})
&=&
\int_{s}^{t}(\int_{\{w=r\}}\frac{1}{|\nabla w|''})\,dr\\
&\geq&
\int_{s}^{t}h(r)u(r)v(r)\,dr.
\end{eqnarray*}
Given $v>0$ and $s>0$, there exists $t$ such that $\vol''(\{s<w<t\})=v$. Then $\vol''(\partial\{s<w<t\})\leq \vol''(\{w=s\})+\vol''(\{w=t\})= uv(s)+uv(t)$ tends to 0 as $s$ tends to $+\infty$. This shows that $I_{(M,g'')}\equiv 0$. Since $\displaystyle\int_{1}^{+\infty}h(t)\,dt=+\infty$, $(M,g'')$ is complete.

If $M$ is real analytic, the construction is the same, in the real analytic category.\qed

\section{The calibration argument}
\label{cali}
\medskip

The following lemma goes back to S.T. Yau, \cite{Yau}.

\begin{lem}
\label{ball}
Let $M$ be a complete $n$-dimensional Riemannian manifold. Let $x\in M$ have injectivity radius $\geq R$ and let $B=B(x,R)$ be the geodesic ball. Assume that on $B$, the sectional curvature is less than $-\rho^2$. Then $I_B (v)\geq (n-1)\rho v$.
\end{lem}

\proof
On $B$, let $r$ denote the distance to $x$, and $\omega=\iota_{\nabla r}vol$. Then $|\omega|\leq 1$ and $d\omega=h\,vol$ where $h$, the mean curvature of geodesic spheres, is controlled by curvature, $h\geq (n-1)\rho$. This form can be used as a ``calibration'' : given a domain $D\subset B$,
\begin{eqnarray*}
(n-1)\rho \vol(D)\leq \int_{D}h\,vol=\int_{\partial D}\omega\leq \lng(\partial D ).\qed
\end{eqnarray*} 

\medskip

In this section, we turn Yau's observation into a systematic tool for computing isoperimetric profiles. Note that, in order to prove isoperimetric inequalities, F. H\'elein has made a slightly different use of calibrations, \cite{Helein}.

\begin{prop}
\label{calib}
Let $M$ be an oriented Riemannian manifold with volume element $vol_{M}$, let $\omega$ be an $n-1$-form on $M$ such that 
\begin{displaymath}
|\omega|\leq 1, \quad\textrm{ and }\quad d\omega=c\,vol_{M}
\end{displaymath}
for some constant $c$.
Let $V\subset M$ be a submanifold with compact boundary and finite volume $v$. Assume that
\begin{itemize}
  \item $\omega$ calibrates $\partial V$, i.e.
\begin{displaymath}
\omega_{|\partial V}=vol_{\partial V}.
\end{displaymath}
  \item there exist compact domains $V_{\epsilon}$ such that
\begin{displaymath}
vol(V\setminus V_{\epsilon})\to 0,\quad vol(\partial V \triangle \partial V_{\epsilon})\to 0
\end{displaymath}
as $\epsilon$ tends to 0.
\end{itemize}
Then $I_{M}(v)=c\,v$.
\end{prop}

\proof

Let $D$ be an arbitrary domain. Then
\begin{displaymath}
c\,vol(D)=\int_{D}c\,vol=\int_{D}d\omega=\int_{\partial D}\omega\leq vol(\partial D)
\end{displaymath}
since $|\omega|\leq 1$. This shows that
\begin{displaymath}
I_{M}(v)\geq c\,v.
\end{displaymath}

Equality holds asymptoticly for sets of the form $V_{\epsilon}$ as $\epsilon$ tends to 0. Indeed, since $\omega$ calibrates $\partial V$,
\begin{displaymath}
\int_{\partial V}\omega = vol(\partial V),
\end{displaymath}
whereas 
\begin{displaymath}
\lim_{\epsilon\to 0}\frac{vol(\partial V_{\epsilon})}{c\,vol(V_{\epsilon})}=1,
\end{displaymath}
showing that 
\begin{displaymath}
I_{M}(v) = (n-1)v.\qed
\end{displaymath}

\begin{cor}
\label{cusp}
Let $M$ be a complete constant curvature $n$-dimensional cusp. Then 
\begin{displaymath}
I_{M}(v)=(n-1)v.
\end{displaymath}
\end{cor}

\proof
Let $\xi$ denote the Busemann vectorfield attached to $\zeta$ (i.e. $\xi(x)$ is the initial speed of the unit speed geodesic starting at $x$ and converging to $\zeta$). Then $\xi$ descends to a well defined unit vectorfield on $M$. Consider the $n-1$-form
\begin{displaymath}
\omega=\iota_{\xi}vol_{M}.
\end{displaymath}
Then 
\begin{displaymath}
|\omega|=|\xi|=1, \quad\textrm{ and }\quad d\omega=(n-1)vol.
\end{displaymath}

Horospheres centered at $\zeta$ descend to hypersurfaces orthogonal to $\xi$ in $M$. Each horoball centered at $\zeta$ descends to an open set $B_{v}$ of finite volume $v$, and $v$ takes every value in $(0,+\infty)$. With the orientation induced from horoballs, horospheres are calibrated by $\omega$. Proposition \ref{calib} applies to $B_{v,\epsilon}=B_{v}\setminus B_{\epsilon}$.\qed

\section{Constructing a singular metric from a Morse function}
\label{singular}

\subsection{Standing assumptions on the level function}
\label{assumptions}
In this section, $M$ will denote an orientable 2-manifold equipped with a smooth function $f$ satisfying the following.
\begin{enumerate}
  \item For all $a\leq b\in \mathbf{R}$, $f^{-1}[a,b]$ is compact.
  \item All critical points of $f$ are nondegenerate of index 1.
  \item The values of $f$ at distinct critical points are distinct integers.
  \item Critical values of $f$ form an interval of $\Z$ containing 0 and 1.
\end{enumerate}

In the terminology of the introduction, $f$ is a level function on $M$.

\subsection{The level graph of $f$}
\label{defgraph}

Let $LG$ be the set of equivalence classes of the equivalence relation
\begin{displaymath}
x\mathcal{R}y \Leftrightarrow x \textrm{ is connected to }y\textrm{ in a level set of }f.
\end{displaymath}

\begin{figure}
\begin{center}
\includegraphics[width=2in]{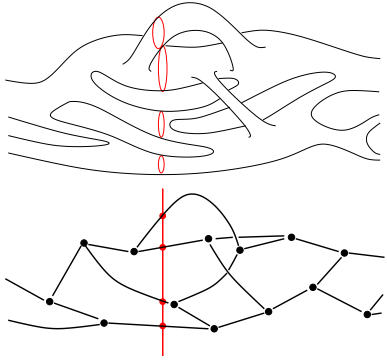}
\caption{The level graph of a function}
\label{graph}
\end{center}
\end{figure}

Under assumptions 1 to 4, $LG$ is a trivalent graph. Indeed, $f$ descends to a map $\f:LG\to \mathbf{R}$ which is a covering map over the complement of critical values. The pull-back of a neighborhood of a critical value consists of finitely many arcs on which $\f$ is a homeomorphism, and a Y shaped set (3 arcs joined at a the class of a critical point), see figure \ref{graph}.
Edges are oriented by the map $\f$.

\subsection{Weights}
\label{defweights}

We are about to put a metric on $M$ such that if $a$ is a regular value of $f$, then $f^{-1}(a)$ is a disjoint union of constant geodesic curvature curves. Each such curve corresponds to an edge of $LG$. The total length of $f^{-1}(a)$ is shared between its connected components in proportion of numbers called the \emph{weights} of the corresponding edges. In this subsection, we describe how these weights are chosen.

\begin{defi}
Weights are given recursively to edges of $LG$ in the following way. Assume that $\f^{-1}(1/2)$ consists of $N$ points. Each of the edges on which $\f$ vanishes is weighted $1/N$. Then as one moves away from $\f^{-1}(1/2)$, following oriented edges, as one encounters a vertex $v$, 
\begin{itemize}
  \item either $v$ has one incoming edge and two outgoing edges; the weight of the incoming edge has been previously defined; then both outgoing edges are weighted one half of the weight of the incoming edge.
  \item or $v$ has two incoming edges and one outgoing edge; the weight of the incoming edges have been previously defined; then the outgoing edge is weighted the sum of the weights of the incoming edges.
\end{itemize}
\begin{figure}
\begin{center}
\includegraphics[width=4in]{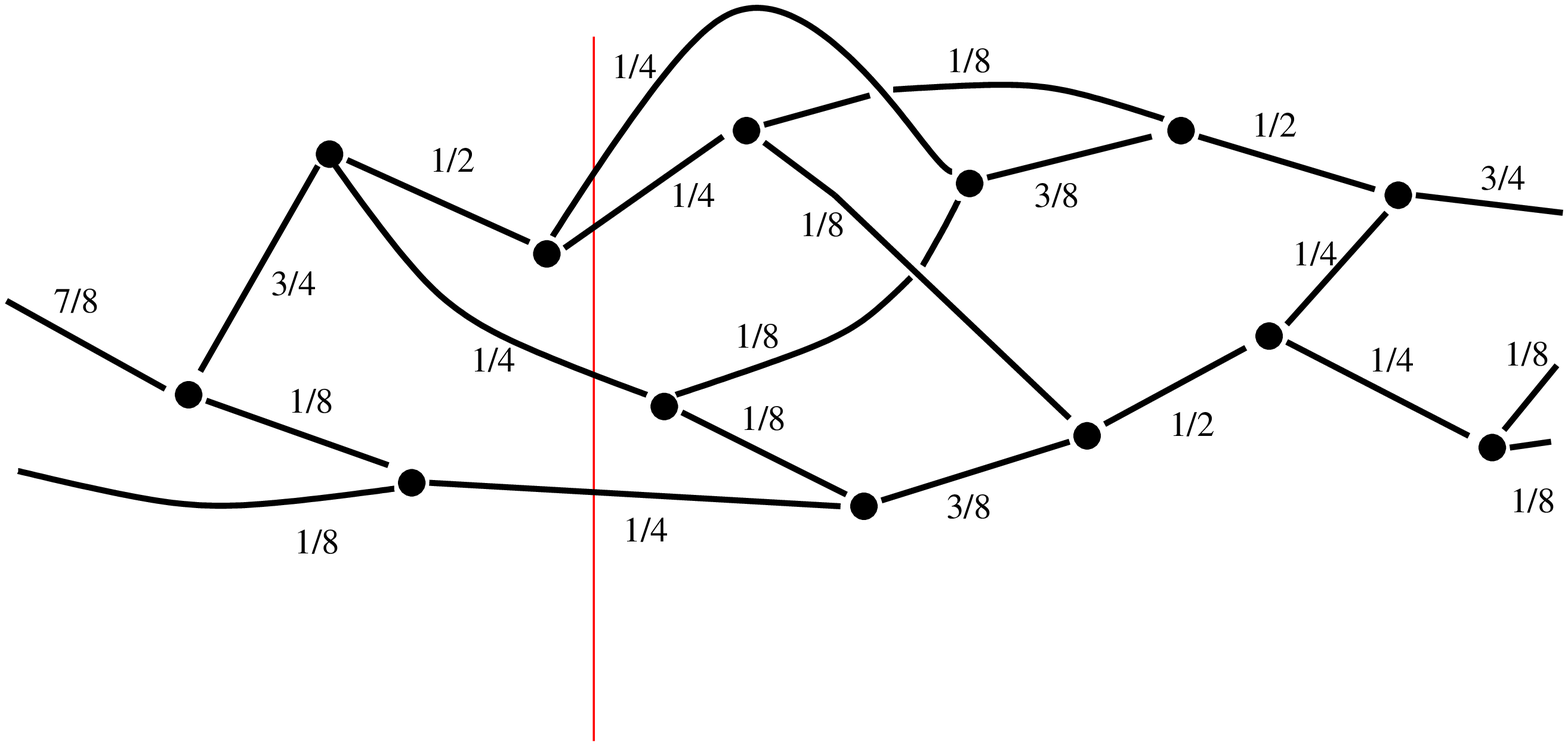}
\caption{Weights}
\label{weights}
\end{center}
\end{figure}
Symmetricly, when proceeding from $\f^{-1}(1/2)$ backwards along oriented edges, the opposite rule is applied : 2 incoming edges share half the weight of 1 outgoing edge, 1 incoming edge adds up the weights of 2 outgoing edges, see figure \ref{weights}.
\end{defi}

Note that for each $t\in\mathbf{R}$, the sum of the weights of the edges on which $\f$ takes the value $t$ is equal to 1. 

\begin{lem}
\label{weight>}
For every oriented edge $e$ of $LG$ originating at $p'=\alpha(e)$ and ending at $p''=\omega(e)$,
\begin{eqnarray*}
weight(e)\geq \frac{1}{N}2^{-|f(p'')|-1}.
\end{eqnarray*}
\end{lem}

\proof
When traversing a vertex of $LG$, weights decrease at worst by a factor of 2. Furthermore, travelling along $LG$ from $e$ to an edge on which $f$ takes the value $1/2$, at most $|f(p'')|+1$ vertices are encountered.\qed

\subsection{Renormalization of the Morse function}
\label{renorm}

Pick an integer $\nu\geq 1$ such that $N\leq 2^{\nu}$.

\begin{defi}
Let $u=\phi\circ f$ where $\phi$ is a diffeomorphism of $\mathbf{R}$ onto $(0,+\infty)$ mapping an integer $n\in\Z$ to $\displaystyle 16^{n(|n|+\nu)}$.
\end{defi}

\begin{lem}
\label{value>}
For every oriented edge $e$ of $LG$ originating at $p'=\alpha(e)$ and ending at $p''=\omega(e)$, and every vertex $p$ such that $u(p)<u(p'')$,
\begin{eqnarray*}
weight(e)u(p'')\geq 8u(p).
\end{eqnarray*}
\end{lem}

\proof
One easily checks that for all $n\in\Z$, $\phi(n)\geq 2^{|n|+4+\nu}\phi(n-1)$. If $f(p'')=n$, then $f(p)\leq n-1$, and Lemma \ref{weight>} gives 
\begin{eqnarray*}
weight(e)&\geq& \frac{1}{N}2^{-|n|-1}\\
&\geq& 2^{-|n|-1-\nu}\\
&\geq&8\frac{\phi(n-1)}{\phi(n)}\\
&\geq&8\frac{u(p)}{u(p'')}.\qed
\end{eqnarray*}

\subsection{Constant curvature annuli}
\label{annuli}

Our model for a calibrated triple $(M,f,\omega)$ is the hyperbolic plane, $f$ a Busemann function and $\omega=*df$ the unit 1-form calibrating a family of concentric horoballs. Denoting $w=e^{-f}$ we have
\begin{displaymath}
\omega=\iota_{w^{-1}\nabla w}vol,
\end{displaymath}
\begin{eqnarray*}
|\nabla w|=w \quad \textrm{and}\quad \Delta w=-2w.
\end{eqnarray*}
In coordinates, $\displaystyle w=\frac{1}{\Im m(z)}=\frac{1}{y}$ when the hyperbolic plane is represented as the upper half plane $\{\Im m(z)>0\}$ equipped with the metric
\begin{displaymath}
\frac{dx^{2}+dy^{2}}{y^{2}}=w^{2}dx^{2}+
\frac{dw^{2}}{w^{2}}.
\end{displaymath}
When divided by a translation in the $x$ coordinate, this simply connected model gives rise to complete constant curvature cusps
\begin{displaymath}
w^{2}dx^{2}+\frac{dw^{2}}{w^{2}},\quad x\in\mathbf{ R}/\tau\mathbf{ Z},\quad w\in(0,+\infty).
\end{displaymath}
and then to annuli
\begin{displaymath}
A_{\tau,c,c'}=(\mathbf{ R}/\tau\mathbf{ Z}\times[c,c'],w^{2}dx^{2}+\frac{dw^{2}}{w^{2}}).
\end{displaymath}
The annulus $A_{\tau,c,c'}$ has constant curvature $-1$. Its boundary components have lengths $c\tau$ and $c'\tau$. Note that $A_{\tau,c,c'}$ is isometric to $A_{c\tau,1,c'/c}$. The second parameter is only used to adjust the interval of variation of the $w$ coordinate. Note that $area(A_{\tau,c,c'})=\tau(c'-c)$, and that the \emph{height} of $A_{\tau,c,c'}$, i.e. the distance between the boundary components, is equal to $\log(c'/c)$. 

\begin{defi}
\label{defanticusp}
We shall call $A_{1,0,\delta}$ (resp. $A_{1,\delta,+\infty}$) a constant curvature $-1$ \emph{cusp} (resp. \emph{anticusp}) with boundary a horocycle of length $\delta$.
\end{defi}

\subsection{Piecing annuli together}
\label{piec}

For each edge $e$ of the level graph $LG$, we use a constant curvature annulus $A(e)=A_{\tau,c,c'}$ where $c=value(\alpha(e))$, $c'=value(\omega(e))$, $\tau=weight(e)$. 

Note that
\begin{displaymath}
area(A(e))=weight(e)(u(\omega(e))-u(\alpha(e))).
\end{displaymath}
The length of the left hand side (resp. right hand side) component of $\partial A(e)$ is 
\begin{displaymath}
length(\partial_{left}A(e))=weight(e)u(\alpha(e)),\quad length(\partial_{right}A(e))=weight(e)u(\omega(e)).
\end{displaymath}
The height of $A(e)$ is 
\begin{displaymath}
\log u(\omega(e))-\log u(\alpha(e)).
\end{displaymath}
In particular, the height is always $\geq \log 2$.

\begin{figure}
\begin{center}
\includegraphics[width=3in]{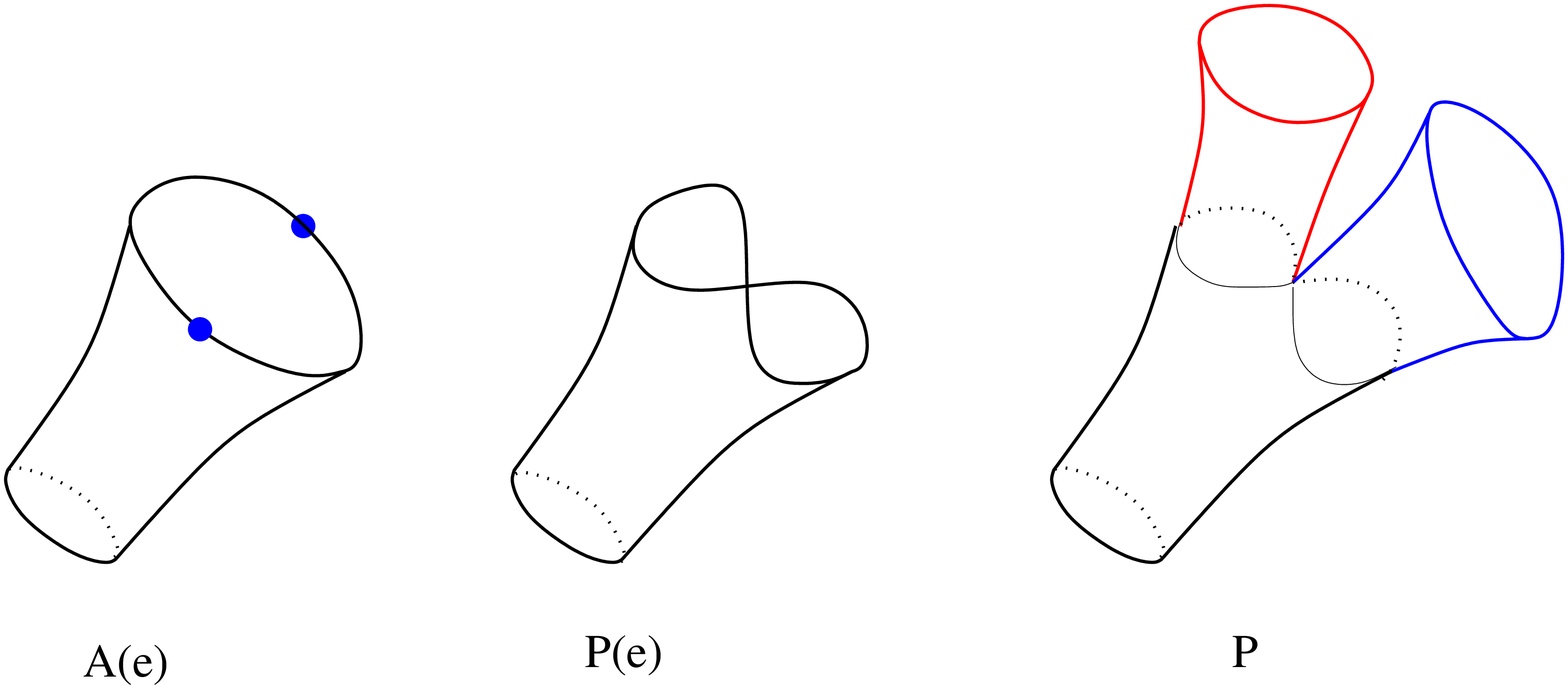}
\caption{Gluing together singular pairs of pants}
\label{P}
\end{center}
\end{figure}

If $\alpha(e)$ (resp. $\omega(e)$) has two incoming (resp. outgoing) edges, one chooses two points on the left hand (resp. right hand) boundary component of $A(e)$ which separate it into two intervals of lengths proportional to the weights of incoming (resp. outgoing) edges. Then the two marked points are identified. The resulting surface $P(e)$ is a singular pair of pants whose boundary consists of three circles, see figure \ref{P}.

By construction, if $e$, $e'$ are incoming and $e''$ is the outgoing edge of some vertex $v$, the lengths of the right hand boundary components of $P(e)$ and $P(e')$ fit with the lengths of the left hand boundary components of $P(e'')$. Therefore one can choose orientation preserving isometries between these circles and glue them together. This applies as well to vertices with one incoming and two outgoing edges.

Morse theory shows that the resulting space
\begin{displaymath}
P=\coprod_{e}P(e)/\sim
\end{displaymath} 
is homeomorphic to $M$, in such a way that $u$ is mapped to $w$, the projection onto the second factor, on each annulus $A(e)$.

\subsection{The nonsingular part of $P$}
\label{nons}

The space $P$ constructed above can be viewed as $M$ equipped with a Riemannian metric $g$ which, away from critical points of $u$, is locally isometric to hyperbolic plane with curvature $-1$. More precisely, for a point $p\in P$, let $\ell(p)$ be the smallest length of a cycle in $\{u=u(p)\}$ containing $p$, let $r(p)$ be the min of $\log 2$ and of the injectivity radius of the bi-infinite constant curvature cusp along the horocycle of length equal to $\ell(p)$. Then the ball $B(p,r(p))$ 
\begin{itemize}
  \item either is isometric to a ball of radius $r(p)$ in the hyperbolic plane,
  \item or contains a critical point.
\end{itemize} 

\begin{prop}
\label{isopsing}
Let $M'$ denote the smooth part of $P$. Then $M'$ is an incomplete surface of constant curvature $-1$, diffeomorphic to $M$ with the critical points of $u$ removed. Its isoperimetric profile is 
\begin{displaymath}
I_{M'}(v)=v.
\end{displaymath}
\end{prop}

\proof
Let $B_{t}=\{u\leq t\}$. Then 
\begin{displaymath}
area(B_{t}\setminus B_{s})=t-s.
\end{displaymath}
In particular $area(B_{t})=t$ is finite.

Indeed, if the interval $[s,t]$ does not contain critical values, $B_{t}\setminus B_{s}$ is the union of disjoint annuli isometric to
\begin{displaymath}
\{s\leq w \leq t\}\subset A_{weight(e),value(\alpha(e),value(\omega(e)}.
\end{displaymath}
Each annulus has area $weight(e)(t-s)$, the weights of edges crossing a level set of $u$ sum up to 1, therefore the total area is $t-s$. The general case follows by additivity. 

If $D$ is a disk of radius $r$ (small enough) centered at a critical point, then $D$ is an isometric double covering of a hyperbolic disk of radius $r$ (see next paragraph). Thus its area and boundary length tend to 0 as $r$ tends to 0. Each $B_{t}$ can be approximated by domains $B_{t,\epsilon}$ by removing $B_{\epsilon}$ and tiny disks around critical points, whose total area and boundary length is arbitrarily small.

Since $u$ is smooth, the 1-form $\omega=\iota_{u^{-1}\nabla u}vol$ is smooth. Since $u$ coincides with coordinate $w$ on each $P(e)$, $\omega$ has unit norm, satisfies $d\omega=vol$ and calibrates the sublevel sets of $u$.
Therefore proposition \ref{calib} applies, and the isoperimetric profile of $M'$ is linear.\qed

\subsection{Description of singularities}
\label{descr}

We claim that at each critical point $p$ of $u$, $P$ has a conical singularity with angle $4\pi$. More precisely, for $r\leq r(p)$, the ball $B(p,r)$ in $M$ is isometric to the disk $D_{r}$ of radius $\sqrt{\tanh(r)}$ in $\mathbf{C}$ equipped with the metric $g_{0}=16|z|^{2}(1-|z|^{4})^{-2}|dz|^{2}$ in such a way that $u=u(p)\frac{|z^{2}-i|^{2}}{1-|z|^{4}}$. 

Indeed, the map $z\mapsto z'=z^{2}$ is a 2-sheeted Riemannian covering of $D_{r}\setminus\{0\}$ onto the punctured disk of radius $\tanh(r)$ equipped with the hyperbolic metric $4|dz'|^{2}/(1-|z|^{2})^{2}$. We use the isometry
\begin{displaymath}
z'\mapsto z''=\frac{1-iz'}{z'-i}
\end{displaymath}
between the Poincar\'e disk and the upper half plane equipped with $|dz''|^{2}/(\Im m(z))^{2}$, and pull back the function $w''=1/\Im m(z'')$. Taking a constant multiple of it, we get a function
\begin{displaymath}
w_{0}=u(p)\frac{|z^{2}-i|^{2}}{1-|z|^{4}}
\end{displaymath}
on $D_{0}$. The level set $\{w_{0}=u(p)\}$ splits $D_{0}$ into four regions, isometric in pairs, and isometric to the intersection of a horoball (resp. the complement of a horoball) and a hyperbolic ball of radius $1/4$ centered on the boundary of the horoball. In each region, $\log(w_{0}/u(p))$ is the signed distance to the boundary horocycle. Therefore these regions fit together exactly into a ball of radius $r$ centered at a singular point of $M$. And $w_{0}$ is exactly mapped to $u$.

\section{Smoothing of singularities}
\label{smoothing}

\begin{lem}
\label{lissage}
Let $D$ be a $n$-dimensional disk, $g_{0}$ a smooth Riemannian metric and $\omega_{0}$ a smooth $n-1$-form on $D$ such that 
\begin{displaymath}
|\omega_{0}|_{0}=1\quad\textrm{and}\quad d\omega_{0}=vol_{0}
\end{displaymath}
in a neighborhood of the boundary. Assume that
\begin{displaymath}
\int_{D}d\omega_{0}>0.
\end{displaymath}
Then there exists a smooth Riemannian metric $g$ and a smooth $n-1$-form  $\omega$ on $D$ such that 
\begin{enumerate}
  \item $g=g_{0}$ and $\omega=\omega_{0}$ in a neighborhood of $\partial D$, 
  \item $|\omega|\leq 1$ and $d\omega=vol$ in $D$,
  \item $vol(D,g)=\int_{D}d\omega_{0}\leq vol_{0}(\partial D)$.
\end{enumerate}
\end{lem}

\proof
First change $g_{0}$ into $g_{1}$ away from the boundary, in order that
\begin{displaymath}
vol(D,g_{1})=\int_{D}d\omega_{0}.
\end{displaymath}
This is possible since $\displaystyle \int_{D}d\omega_{0}>0$ by assumption.

Let $\eta$ be a compactly supported $n$-form on $D$. According to the Poincar\'e Lemma, there exists a compactly supported $n-1$-form $\beta$ on $D$ such that $d\beta=\eta$ if and only if $\int_{D}\eta=0$. Therefore there exists a compactly supported $n-1$-form $\beta$ such that $\omega=\omega_{0}+\beta$ satisfies
\begin{displaymath}
d\omega=vol_{1}.
\end{displaymath}

Then we change $g_{1}$ into $g$ while keeping the same volume element. At points where $\omega$ does not vanish, one can write
\begin{displaymath}
g_{1}=(\frac{*_{1}\omega}{|\omega|_{1}})^{2}+g_{1}^{\perp}
\end{displaymath}
where $g_{1}^{\perp}$ is non negative with kernel the vectorfield $\xi$ dual to $\omega$. 

Let
\begin{displaymath}
g=(\frac{*_{1}\omega}{\chi(|\omega|_{1})^{-1}|\omega|_{1}})^{2}+\chi(|\omega|_{1})^{-2/(n-1)}{g_{1}^{\perp}},
\end{displaymath}
where $\chi:\mathbf{R}_{+}\to(0,+\infty)$ is a smooth function such that
\begin{itemize}
  \item $\chi(t)\geq t$ ;
  \item $\chi(t)=t$ for $t>2/3$ ;
  \item $\chi(t)=1$ for $t<1/3$.
\end{itemize}
With respect to the metric $g$, the 1-form $\displaystyle \frac{*_{1}\omega}{\chi(|\omega|_{1})^{-1}|\omega|_{1}}$ has unit norm, so that
\begin{displaymath}
|\omega|=|*\omega|=\frac{|\omega|_{1}}{\chi(|\omega|_{1})}\leq 1.
\end{displaymath}
In a neighborhood of the boundary, $g_{1}=g_{0}$ and $\omega=\omega_{0}$, thus $|\omega|_{1}=1$ and $g=g_{1}$. At points where $|\omega|_{1}<1/3$, $g=g_{1}$, thus $g$ extends smoothly to all of $D$, and the inequality $|\omega|\leq 1$ trivially extends to points where $\omega$ vanishes.

Since $vol_{g}=vol_{1}$, $d\omega=vol_{g}$ and \begin{eqnarray*}
area(D,g)
&=&area(D,g_{1})\\
&=&\int_{D}d\omega_{0}\\
&=&\int_{\partial D}\omega_{0}\\
&\leq&vol_{0}(\partial D).\qed
\end{eqnarray*}

\begin{prop}
\label{smooth}
Let $E$ be the union of the intervals $(c/2,2c)$ where $c$ is a critical value of $u$. There exists a smooth Riemannian metric and a smooth 1-form $\omega$ on $M$ with the following properties.
\begin{enumerate}
  \item $|\omega|\leq 1$.
  \item $d\omega=vol$.
  \item $\omega$ calibrates the sublevel sets $\{u\leq t\}$ for all $t\notin E$.
\end{enumerate}
Furthermore, one can arrange that each superlevel set $\{u\geq s\}$ have bounded geometry.
\end{prop}

\proof
This local construction can be performed in the model $(D_{r},g_{0})$. In order to apply Lemma \ref{lissage} to smooth versions of $g_{0}$ and $\omega_{0}$, we must check that they satisfy
\begin{displaymath}
I=\int_{D_{r}}d\omega_{0}>0.
\end{displaymath}
Since
\begin{displaymath}
\int_{D_{0}}d\omega_{0}
=\int_{\partial D_{0}}\omega_{0},
\end{displaymath} 
$I$ only depends on the boundary data, and we can compute it from the singular data. Using the double covering $z\mapsto z^{2}$ of $D_{r}$ onto a ball of hyperbolic radius $\tanh(r)$, we see that $I$ is twice the flux of a Busemann vectorfield along the boundary of a hyperbolic ball, i.e. twice the area of the hyperbolic ball. In particular, $I>0$.

Therefore, Lemma \ref{lissage} can be applied near each critical point $p$. Let $\ell(p)$ denote the minimum length of a boundary component of an annulus that contains $p$. Then $\ell(p)$ increases with $u(p)$. Indeed, when moving from a critical value to the next, the total length of the critical level is multiplied at least by 16, and weights are at most divided by 2. The
injectivity radius of the bi-infinite curvature $-1$ cusp at a point on the horocycle of length $\ell(p)$ also increases. Therefore, the radius $r(p)$ of the ball centered at a critical point $p$ to which Lemma \ref{lissage} can be applied increases with $u(p)$. On the ball $B(p,r(p))$, $u$ takes values in $(e^{-r(p)}u(p),e^{r(p)}u(p))$. Since $E$ contains the interval $(u(p)/2,2u(p))$, one can perform exactly the same smoothing on $B(p,\log 2)$ for all $p$, provided $u(p)$ is large enough. This gives bounded geometry.\qed

\begin{rem}
If $c$ is a critical value of $u$, then $D=\{u\leq c\}$ is \emph{not} a solution of the isoperimetric problem.
\end{rem}
Indeed, solutions have a smooth boundary. Compare \cite{Calabi-Cao} page 541 or \cite{Pitts} page 37.

\section{Construction of special Morse functions}
\label{morse}

\begin{prop}
\label{special}
Let $M$ be a connected orientable 2-manifold. Assume that one of the following properties holds.
\begin{itemize}
  \item $M$ has at least 4 ends.
  \item $M$ has 3 ends and at least one of the ends has infinite genus.
  \item $M$ has 2 ends and both have infinite genus. 
\end{itemize}
Then there exists on $M$ a function $f$ such that
\begin{enumerate}
  \item For all $a\leq b\in\mathbf{R}$, $f^{-1}[a,b]$ is compact.
  \item All critical points of $f$ are nondegenerate of index 1.
  \item The values of $f$ at distinct critical points are distinct.
  \item Critical values of $f$ form a discrete set which does not contain 0.
  \item All level sets of $f$ are disconnected.
\end{enumerate}
\end{prop}

\begin{lem}
\label{43}
Let $M$ be a connected orientable compact 2-manifold with boundary split as $\partial M=\partial_{-}M\cup \partial_{+}M$ where both pieces are nonempty unions of connected components. There exists a Morse function $f$ on $M$ which is constant equal to 1 (resp. $-1$) on $\partial_{+}M$ (resp. $\partial_{-}M$), takes values in $[-1,1]$, has distinct critical values and no critical points of index 0 or 2.

If furthermore both $\partial_{+}M$ and $\partial_{-}M$ are disconnected, then one can arrange that all levels of $f$ be disconnected. 
\end{lem}

\begin{figure}
\begin{center}
\includegraphics[width=2in]{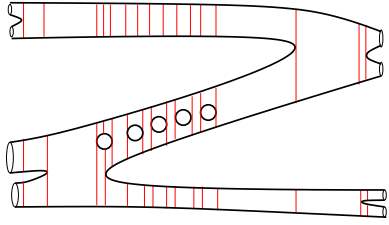}
\caption{A Morse function with disconnected levels on a compact surface with boundary}
\label{4}
\end{center}
\end{figure}

\proof
Compact orientable surfaces with marked boundary are classified by their genus and the number of boundary components of each sign. Indeed, such surfaces are obtained by removing disks from closed orientable surfaces. Closed orientable surfaces are classified by their genus and the diffeomorphism group of a surface is transitive on finite collections of disks. Therefore is suffices to give an example for each value of the triple (genus, number of $+$ boundary components, number of $-$ boundary components). Figure \ref{4} shows such a surface together with the critical levels of a suitable Morse function. 
If $\partial_{+}M$ or $\partial_{-}M$ is connected, the Z-shaped core must be replaced with a $>$ or --- -shaped surface, in which case the disconnectedness of levels is lost.\qed

\begin{lem}
\label{index0}
Let $M$ be a connected orientable non compact 2-manifold with non empty compact boundary. There exists a proper Morse function on $M$ which is constant along the boundary and satisfies 1 to 4.
\end{lem}

\proof
Start with an arbitrary proper Morse function $f$ which is constant on the boundary. Choose an increasing sequence $t_{j}$ of noncritical values. On each compact submanifold $\{t_{j}\leq f\leq t_{j+1}\}$, apply Lemma \ref{43} to improve $f$. The resulting function has all the required properties.\qed

\subsection{Proof of Proposition \ref{special}}

Assume $M$ has at least $k$ ends. Then for every sufficiently large compact subset $K\subset M$, $M\setminus K$ has at least $k$ non compact connected components. Choose such a $K$ with smooth boundary. For each connected component $C$ of $M\setminus K$, $\bar{C}\cap K$ is open and closed in $K$. Therefore there are finitely many components $C$. Add the compact ones to $K$. Using Lemma \ref{index0}, choose on each non compact component $C$ a proper Morse function which is constant equal to 1 on $\partial C$ and satisfies 1 to 4. Change the sign of the Morse function on two of the non compact connected components. Apply Lemma \ref{43} to get a Morse function on all of $M$.

Case $k=4$. For $t\in[-1,1]$, the level $\{f=t\}$ is disconnected as shown on the picture. For $t>1$ or $t<1$, this follows from the fact that $\{f=t\}$ intersects at least two connected components of $M\setminus K$.

\begin{figure}
\begin{center}
\includegraphics[width=2in]{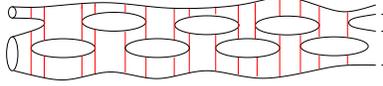}
\caption{Infinite genus end}
\label{infinitegenusend}
\end{center}
\end{figure}

Case $k=3$. We can assume that the end where $u\geq 1$ has infinite genus. A noncompact end of infinite genus is, up to a compact change, diffeomorphic to the surface with boundary shown on figure \ref{infinitegenusend} (\cite{Richards}, \cite{Kerekjarto}). Therefore, up to enlarging $K$, we can assume that the connected component of $M\setminus K$ where $f\geq 1$ has a boundary consisting of exactly 2 curves.
Figure \ref{infinitegenusend} also shows the critical levels of a Morse function on the end satisfying 1 to 5. We use it to modify $f$. Thus if $t\geq 1$, the level $\{f=t\}$ is disconnected. For $t\in[-1,1]$, this follows from Lemma \ref{43}. For $t\leq 1$, this follows from the fact that $\{f=t\}$ intersects two connected components of $M\setminus K$.

Case $k=2$. An orientable surface with two ends of infinite genus is the bi-infinite version of figure \ref{infinitegenusend} (\cite{Richards}, \cite{Kerekjarto}) and therefore has an obvious Morse function satisfying 1 to 5.\qed

\section{Proof of Theorem \ref{main}}
\label{proof}

Apply Proposition \ref{special} in order to obtain a Morse function on $M$ having properties 1 to 5. Apply Proposition \ref{smooth} to produce a smooth metric, a smooth function $u$ and a smooth calibration $\omega$ on $M$ which calibrates the boundary of most sublevel sets 
\begin{displaymath}
B_{t}=\{u\leq t\}
\end{displaymath}
of $u$, but unfortunately not all of them. This works for levels which are not in an open neighborhood $E$ of the set of critical values. $E$ is a discrete union of intervals of the form $(c_{-},c_{+})$, each of them containing exactly one critical value $c$, $c_- =c/2$, $c_+ =2c$. 

For each $t\in (c_{-},c_{+})$, we must find a submanifold with boundary which is calibrated by $\omega$ and has area $t$. 
 
\begin{figure}
\begin{center}
\includegraphics[width=3in]{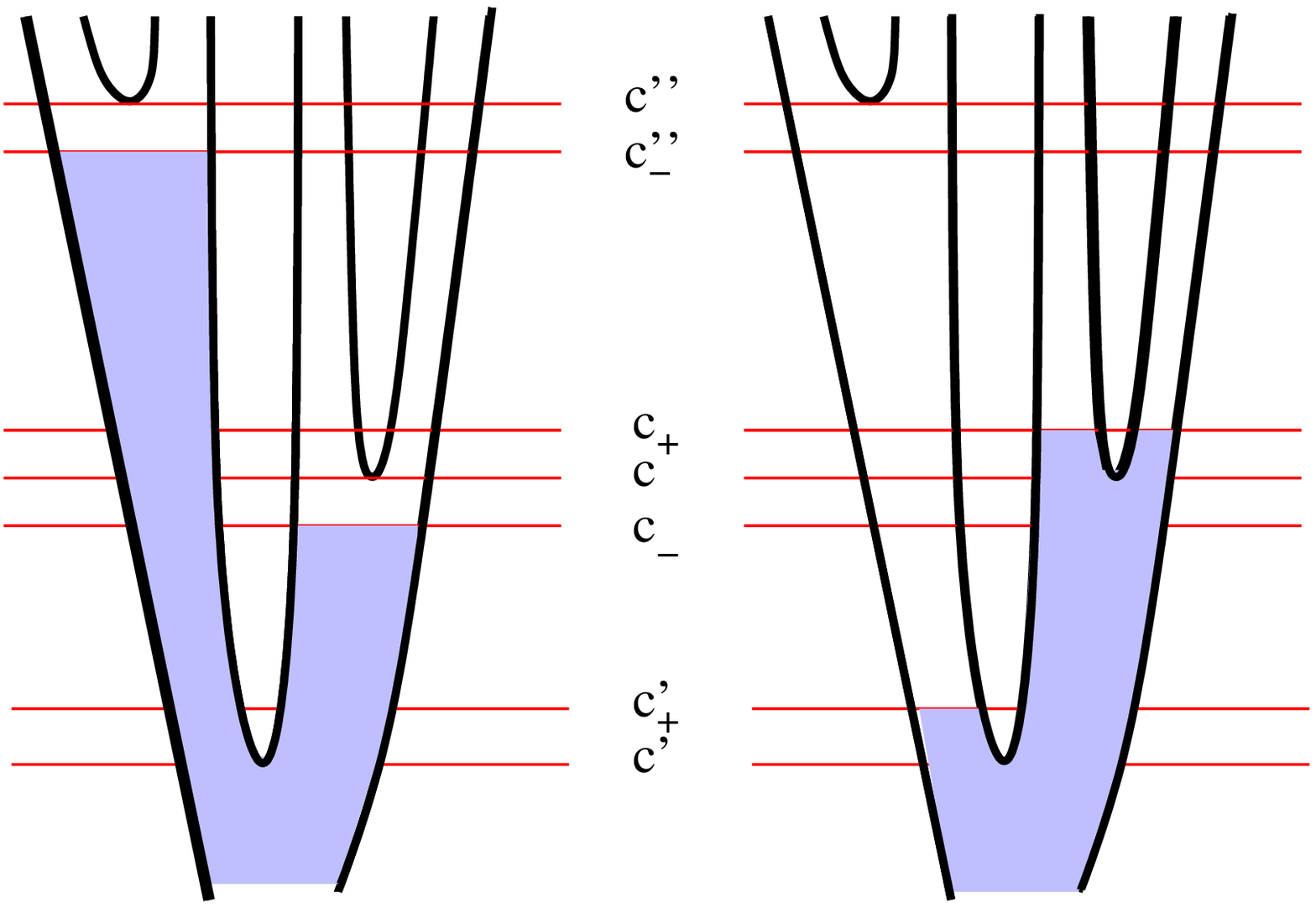}
\caption{$C_{e,c''_{-}}$ \hskip2cm $D_{e,c'_{+}}$}
\label{vase}
\end{center}
\end{figure}

We use the fact that the level set $L=\{u=c\}$ is disconnected. Let $L'$ denote a connected component of $L$ which does not contain the critical point, and $L''=L\setminus L'$. $L'$ is an element of an edge $e$ of the level graph of $u$. Let $c'=value(\alpha(e))$ and $c''=value(\omega(e))$. 
For $s\in (c'_{+},c''_{-})$, let
\begin{displaymath}
C_{e,s}=B_{c_{-}}\cup (P(e)\cap\{w\leq s\})
\quad\hbox{ and }\quad
D_{e,s}=B_{c_{+}}\setminus (P(e)\cap\{w> s\}),
\end{displaymath}
see figure \ref{vase}.
Then $\partial C_{e,s}$ and $\partial D_{e,s}$ are calibrated by $\omega$. 
Using the fact that $c''\geq 16 c'$ and Lemma \ref{value>}, one gets
\begin{eqnarray*}
area(C_{e,c''_{-}})-
area(D_{e,c'_{+}})
&\geq&
area(P(e)\cap(B_{c''_{-}}\setminus B_{c'_{+}}))-area(B_{c_{+}}\setminus B_{c_{-}})\\
&=&
weight(e)(c''_{-}-c'_{+})-(c_{+}-c_{-})\\
&\geq &
weight(e)(\frac{1}{2}-\frac{1}{8})c''-(2-\frac{1}{2})c\\
&\geq&
\frac{1}{4}weight(e)c''-2c\\
&\geq&0,
\end{eqnarray*}

Therefore, as $s$ varies, $area(C_{e,s})$ and $area(D_{e,s})$ together  achieve all values in $(c_{-},c_{+})$ (this is the pipe clearing trick). Then Proposition \ref{calib} applies, and $I_{M}(a)=a$ for all $a>0$.\qed

\section{Surfaces of revolution}
\label{revolution}

This is a preparation for the next section, where caps of revolution with explicit isoperimetric profiles will be needed.

In a surface of revolution, the isoperimetric profile is not always achieved by disks of revolution. It is the case provided the curvature is a nonincreasing function of the distance from the pole. This is not sufficient for our purposes. Indeed, a constant curvature cusp cannot be replaced with a surface of revolution with nonincreasing curvature of the same area. We need just a little more flexibility. Fortunately, the property that disks of revolution are extremal is often stable under small perturbations.

\begin{lem}
\label{unique}
Let $A=([a,b]\times \R/2\pi\Z,dr^2 +f(r)^2 d\theta^2 )$, $f>0$, be a smooth annulus of revolution such that for all $r\in[a,b]$, $f'(r)^2-f''(r)f(r)$ is not the square of a nonzero integer. There exist constants $\beta>0$ and $\gamma>0$ such that if $g$ is a $C^2$ function on $[a,b]$ such that
$\n{g-f}_{C^2}<\beta$, $D_0 =\{a\leq r<r_0\}$ and $D=\{a\leq r<\rho(\theta)\}$ are domains such that $\area(D)=\area(D_0 )$ and $\n{\rho -r_0 }_{C^2}<\gamma$ whose boundary has constant geodesic curvature in the metric $dr^2 +g(r)^2 d\theta^2$ , then $D=D_0$.
\end{lem}

\proof
Given $g\in C^2 ([a,b])$ and $\rho\in C^2 (\R/2\pi\Z)$, let $\kappa\in C^0 (\R/2\pi\Z)$ be the geodesic curvature of the curve $\{r=\rho(\theta)\}$,
\begin{eqnarray}
\kappa=\frac{-g(\rho)\rho''+g'(\rho)g(\rho)^2 +2g'(\rho)\rho'^2}{\sqrt{g(\rho)^2 +\rho'^2}^3}.\label{curv}
\end{eqnarray}
The area of $D=\{a\leq r<\rho(\theta)\}$ is equal to 
\begin{eqnarray*}
\area(D)
&=&
\int_{0}^{2\pi}\int_{a}^{\rho(\theta)}g(r)\,dr\,d\theta\\
&=&
\int_{0}^{2\pi}G(\rho(\theta))\,d\theta,
\end{eqnarray*}
where $G(r)=\int_{a}^{r}g(s)\,ds$. Therefore $C^2$ metrics and $C^2$ annular domains of a given area $v$ are conveniently described by pairs of functions $G:[a,b]\to\R$, $R:\R/2\pi\Z\to\R$ such that $G(0)=0$, $G'>0$ and $\int_{0}^{2\pi}R(\theta)\,d\theta=0$. $g$ and $\rho$ are recovered via $g=G'$ and $\rho=G^{-1}\circ (R+v)$. Let $C_0^3 ([a,b])$ denote the space of $C^3$ functions $G$ on $[a,b]$ such that $G(a)=0$. Let $C_0^2(\R/2\pi\Z)$ denote the space of $C^2$ functions $R$ such that $\int_{0}^{2\pi}R(\theta)\,d\theta=0$. Let $U\subset \R\times C_0^3 ([a,b])\times C_v^2(\R/2\pi\Z)$ denote the open subset of triples $(v,G,R)$ such that $v>0$, $G'>0$ and $0<R+v<G(b)$. Let $C_0^0 (\R/2\pi\Z)$ denote the space of continuous functions with vanishing average. Let $\Phi:U\to C_0^0 (\R/2\pi\Z)$ denote the $C^1$ map defined by
\begin{eqnarray*}
\Phi(v,G,R)=\kappa-\frac{1}{2\pi}\int_{0}^{2\pi}\kappa(\theta)\,d\theta,
\end{eqnarray*}
where $\kappa=\kappa(g,\rho)$ is given by formula (\ref{curv}), $g=G'$ and $\rho=G^{-1}\circ (R+v)$.

Constant geodesic curvature annular domains of area $v$ correspond to solutions of $\Phi(v,G,R)=0$. For every $(v,G)$ such that $0<v<G(b)$, $R=0$ is a solution. Let us show that, for $(v,G)$ close to $(v,F)$ where $F(r)=\int_{0}^{r}f(s)\,ds$, this is the only solution. Let us compute the derivative of $\Phi$ with respect to $R$ at $(v,F,0)$. Let $r_0 =F^{-1}(v)$ be the solution of $F(r_0 )=v$. Let $\rho_t$ be a smooth family of functions starting from the constant function $\rho_0 =r_0$, and $\displaystyle \dot{\rho}=\frac{d\rho_t}{dt}_{|t=0}$. Then 
\begin{eqnarray*}
\frac{d\Phi}{dt}_{|t=0}
&=&
-3f^{-4} f'\dot{\rho}(f'f^2 )+f^{-3}(-f\dot{\rho}''+f''f^2 \dot{\rho}+2f'^2 f\dot{\rho})
\end{eqnarray*}
where $f=f(r_0 )$, $f'=f'(r_0 )$, and $f''=f''(r_0 )$. In other words, $\displaystyle \frac{\partial\Phi}{\partial R}(v,F,0):C_0^2 (\R/2\pi\Z)\to C_0^0 (\R/2\pi\Z)$ is equal to the operator 
\begin{eqnarray*}
P:u\mapsto -f^{-2}u''+(f''f^{-1}-f'^2 f^{-2})u.
\end{eqnarray*}
Since the eigenvalues of $u\mapsto -u''$ on functions on $\R/2\pi\Z$ with vanishing average are squares of nonzero integers, the assumption made implies that $P$ is invertible on Sobolev spaces $H_0^2 \to H_0^0$ of functions with vanishing average which are in $L^2$ (resp. have 2 derivatives in $L^2$). Since $H_0^2 \subset C_0^0$, if $Pu$ is continuous, so is $u''$, thus $P^{-1}$ maps $C_0^0$ to $C_0^2$. $P$ is a continuous bijection $C_0^2 (\R/2\pi\Z)\to C_0^0 (\R/2\pi\Z)$, and thus an isomorphism. The implicit function theorem applies : there are neighborhoods $\mathcal{I}_{v}$ of $v$ in $\R$, $\mathcal{V}$ of $F$ in $C_0^3 ([a,b])$ and $\mathcal{W}$ of $0$ in $C_0^2 (\R/2\pi\Z)$ such that for $(v',G)\in\mathcal{I}_{v}\times\mathcal{V}_v$, the equation $\Phi(v',G,R)=0$ admits a unique solution in $\mathcal{W}_v$, this is $R=0$. Covering the interval $[0,F(b)]$ with a finite number of intervals $\mathcal{I}_{v}$, this translates into constants $\beta>0$ and $\gamma>0$ such that if $\n{g-f}_{C^2}<\beta$ and $\n{\rho -r_0}_{C^2}<\gamma$, then the curve $\{r=\rho(\theta)\}$ has constant geodesic curvature if and only if $\rho$ is constant.\qed

\subsection{Optimality of rotationally symmetric disks is often stable}

Although we shall not need the following proposition, we state it for its independant interest.

\begin{prop}
\label{profrev}
Let $S=(\R_+ \times \R/2\pi\Z ,dr^2 +f(r)^2 d\theta^2 )$ be a complete smooth surface of revolution such that
\begin{itemize}
  \item $f$ is nondecreasing and tends to $+\infty$;
  \item the curvature $\displaystyle K(r)=-\frac{f''(r)}{f(r)}$ tends to $-\infty$;
  \item the disks of revolution $\{r<t\}$ are the only extremal domains for the isoperimetric problem in $S$, and for all $t\in[a,b]$, they are strictly stable, i.e. 
\begin{eqnarray*}
\frac{f'(t)^2}{f(t)^2}-\frac{f''(t)}{f(t)}<\frac{1}{f(t)^2}.
\end{eqnarray*}
\end{itemize}
Let $g$ be a $C^2$ function on $\R_+$ such that $g=f$ outside $[a,b]$. If $\n{g-f}_{C^2}$ is small enough, then in the surface of revolution $S'=(\R^2 ,dr^2 +g(r)^2 d\theta^2 )$, the disks of revolution $\{r<t\}$ are the only extremal domains for the isoperimetric problem. In particular, for all $t>0$,
\begin{eqnarray*}
I_{S'}(2\pi \int_{0}^{t}g(s)\,ds)=2\pi g(t).
\end{eqnarray*}
\end{prop}

\proof
The assumptions on $f$ imply that all surfaces under consideration are ultrahyperbolic. Therefore extremal domains exist and form compact sets.

The proof goes by contradiction. Assume there exists a sequence $g_{\ell}$ of functions $C^{2}$-converging to $f$ and a sequence of non rotationally symmetric extremal domains $D_{\ell}$. One can assume that $D_{\ell}$ converges in flat norm to some domain $D$. By lower semi-continuity of the boundary length and continuity of the area, $D$ is extremal in $S$, thus $D$ is rotationally symmetric, i.e. $D=\{r\leq r_0 \}$. 

Let $\chi:\R_+ \to \R_+$ be a smooth function such that $\chi(r_0 )=0$. By flat convergence,
\begin{eqnarray*}
\lim_{\ell\to+\infty} \int_{\partial D_{\ell}}\chi(r)f(r)\,d\theta=\int_{\partial D}\chi(r)f(r)\,d\theta=0.
\end{eqnarray*}
This shows that, for every $\e>0$, outside the $\e$ neighborhood $U_{\e}$ of $\partial D$, 
\begin{eqnarray*}
\lim_{\ell\to+\infty}\int_{\partial D_{\ell} \setminus U_{\e}}f(r)\,d\theta =0.
\end{eqnarray*}
In particular,
\begin{eqnarray*}
\lim_{\ell\to+\infty}\lng(\partial D_{\ell} )
&=&
\lng(\partial D)\\
&=&
\int_{\partial D}f(r)\,d\theta\\
&=&
\lim_{\ell\to+\infty} \int_{\partial D_{\ell} \cap U_{\e}} f(r)\,d\theta\\
&\leq&
\lim_{\ell\to+\infty}\lng(\partial D_{\ell} \cap U_{\e}),
\end{eqnarray*}
which shows that $\displaystyle \lim_{\ell\to+\infty}\lng(\partial D_{\ell} \setminus U_{\e})=0$.

Since the projection $(r,\theta)\to r$ is length decreasing, there exist two sequences $r_{\ell} <r_0 <r'_{\ell}$, tending to $r_0$ such that $\partial D_{\ell}$ does not intersect the parallels $\{r=r_{\ell} \}$ and $\{r=r'_{\ell} \}$. Since
\begin{eqnarray*}
\lim_{\ell\to+\infty}\int_{\partial D_{\ell} \cap \{r_{\ell} <r<r'_{\ell} \}}d\theta =\int_{\partial D}d\theta=2\pi,
\end{eqnarray*}
exactly one component $c_{\ell}$ of $\partial D_{\ell}\cap \{r_{\ell} <r<r'_{\ell} \}$ is homotopic to $\partial D$ in the annulus $\{r_{\ell} <r<r'_{\ell} \}$. Using an arclength parametrization of $c_{\ell}$ and applying Ascoli's theorem, one can assume that $c_{\ell}$ converges uniformly, and the limit must be $\partial D$ parametrized by arclength.

$\partial D_{\ell}$ has constant geodesic curvature $\kappa_{\ell}$. At points where the function $r$ restricted to $c_{\ell}$ achieves its maximum $r_{max}$ or its minimum $r_{min}$, a comparison principle yields
\begin{eqnarray*}
\frac{f'(r_{min})}{f(r_{min})}\leq \kappa_{\ell} \leq
\frac{f'(r_{max})}{f(r_{max})},
\end{eqnarray*}
showing that $\kappa_{\ell}$ tends to $\frac{f'(r_0 )}{f(r_0 )}$.

Since short closed curves must have large geodesic curvature somewhere, $\partial D_{\ell}$ has no short components, thus $\partial D_{\ell} =c_{\ell}$.

Along a long curve with bounded geodesic curvature which makes a large angle with parallels at some point, the function $r$ varies a lot. This shows that the angle of the tangent of $c_{\ell}$ with parallels tends uniformly to 0. In other words, $c_{\ell}$ can be viewed as a graph $\{r=\rho_{\ell} (\theta)\}$ where $\rho_{\ell}$ converges to  $r_0$ in $C^1$. Since the geodesic curvature converges, formula (\ref{curv}) for the geodesic curvature shows that $\rho''_{\ell}$ converge, i.e. $\rho_{\ell}$ converges to  $r_0$ in $C^2$. Lemma \ref{unique} implies that, for $\ell$ large enough, $\rho_{\ell}$ is constant, a contradiction.\qed

\subsection{Finite approximations of a cusp-like surface}

What we need is a variant of Proposition \ref{profrev} where the model surface is non compact.

\begin{prop}
\label{profcusp}
Let $S=(\R\times\R/2\pi\Z,dr^2 +f(r)^2 d\theta^2 )$ denote a complete surface of revolution. Let $S_{\ell}=([-R_{\ell},+\infty)\times\R/2\pi\Z,dr^2 +g_{\ell}(r)^2 d\theta^2 )$ be a sequence of complete smooth surfaces of revolution converging to $S$ in the following sense : $g_{\ell}=f$ on $\R_{+}$, and $g_{\ell}$ converges to $f$, $C^{2}$-uniformly on compact subsets of $\R_{-}$. Assume that
\begin{itemize}
  \item $f(r)=e^{r}$ for $r\leq 0$;
  \item $f'/f$ is nondecreasing;
  \item $f$ is nondecreasing and tends to $+\infty$;
  \item the curvature $\displaystyle K(r)=-\frac{f''(r)}{f(r)}$ tends to $-\infty$;
  \item $g_{\ell}$ is nondecreasing.
\end{itemize}
Then for $\ell$ large enough, the disks of revolution $\{r\leq t\}$ are extremal domains for the isoperimetric problem in $S_{\ell}$. In particular, for all $t>-R_{\ell}$,
\begin{eqnarray*}
I_{S_{\ell}}(2\pi \int_{-R_{\ell}}^{t}g_{\ell}(s)\,ds)=2\pi g_{\ell}(t).
\end{eqnarray*}
\end{prop}

\proof
Since the candidate isoperimetric profile is nondecreasing, one can restrict to domains which are disjoint unions of disks (otherwise, replace a domain with the largest disk spanned by one of its boundary components). Since $S_{\ell}$ is ultrahyperbolic, there exist extremal domains, which again are disjoint unions of disks. Let $D_{\ell}$ be extremal domains with area $v_{\ell}$ in $S_{\ell}$ which converge in flat norm to a possibly non compact submanifold $D\subset S$. Ultrahyperbolicity holds for $S$ on $\{r\geq 0\}$. Therefore $D\cap \{r\geq 0\}$ is compact. $D$ can be approximated by domains of the form $D\cap\{e^{r}\leq \e\}$ with $\e$ tending to $0$. Let $\omega=f(r)\,d\theta$ on $S$. Then $d\omega=u\,vol$ where $u=f'/f$ is nondecreasing. Proposition \ref{rear} applies, showing that $D$ has calibrated boundary, i.e. $D=\{r\leq r_{0}\}$ is a rotationally symmetric non compact annulus.

As in the proof of Proposition \ref{profrev}, most of the length of $\partial D_{\ell}$ concentrates in thin neighborhoods of $\partial D$. Then one component $c_{\ell}$ of $\partial D_{\ell}$ is contained in a thinner and thinner neighborhood of $\partial D$ (being homotopic to $\partial D$ in this neighborhood). One can assume that $c_{\ell}$ converges uniformly, its constant geodesic curvature converges, thus $c_{\ell}$ is a graph which $C^{2}$-converges to $\partial D$. Since $f'/f$ is increasing, $f'^{2}-f''f\leq 0<1$.
Lemma \ref{unique} implies that $c_{\ell}=\{r=r_{\ell}\}$ is a parallel for $\ell$ large enough.

For $\ell$ large, $D_{\ell}$ contains $\{r\leq r_{\ell}\}$ with $r_{\ell}$ close to $r_{0}$. The other components of $D_{\ell}$ are thus contained in $\{r\geq r_{\ell}\}$. By uniform ultrahyperbolicity on the $\{r\geq 0\}$ side, they are contained in fact in a fixed compact set on with the metrics converge $C^{2}$-uniformly. Short closed curves contained in this set cannot have bounded geodesic curvature. We conclude that $D_{\ell}=\{r\leq r_{\ell}\}$.\qed

\subsection{Caps with prescribed isoperimetric profile}

\begin{cor}
\label{prescrcaps}
Let $\delta>0$ be small enough and $k$, $\alpha$ be large enough. There exists a smooth complete surface of revolution $S=S_{\delta,k,\alpha}$ such that
\begin{enumerate}
  \item $S$ is ultrahyperbolic;
  \item $S$ contains a rotationally symmetric annulus of constant curvature $-1$ and area $\alpha+\delta$, whose inner boundary has length $\delta$ and inner disk has area $\delta$;
  \item the curvature at the origin is equal to $k$;
  \item all rotationally symmetric disks in $S$ are extremal;
  \item the square $I_{S}^2$ of the isoperimetric profile of $S$ extends to a smooth function on a neighborhood of $\R_+$;
  \item $I_S$ is nondecreasing and $v\mapsto I_{S}(v)/v$ is nonincreasing;
  \item for $\delta\leq v\leq \alpha$, $I_{S}(v)=v$;
  \item for all $v$, $I_{S}(v)\geq \sqrt{4\pi v -k\,v^2}$.
\end{enumerate}
\end{cor}

\proof
Let $S=(\R_+ \times\R/2\pi\Z ,g=dr^2 +f(r)^2 d\theta^2 )$ be a surface of revolution. Let us first collect necessary conditions on the function $f$. Assume that all disks of revolution are extremal domains in $S$. In other words, the isoperimetric profile $I_S$ is determined by
\begin{eqnarray}
V'(r)=I_S (V(r)).\label{VIS}
\end{eqnarray}
where $V(r)=2\pi \int_{0}^{r}f(s)\,ds$ is the area of the disk of revolution of radius $r$. For $f$ to give rise to a smooth metric on a plane, it is necessary and sufficient that $f$ extends to a smooth odd function on $\R$ such that $f'(0)=1$. Therefore $V$ extends to a smooth even function, i.e. $V(r)=G(r^2 )$ where $G$ is smooth on a neighborhood of $0$. Since $V''(0)=2\pi f'(0)=2\pi$, $V(r)\sim \pi r^2$ and $G'(0)=\pi$, so the inverse map $G^{-1}$ is smooth on a neighborhood of 0. Equation (\ref{VIS}) implies that $I_S (G(r^2 ))^2 =4r^2 G'(r^2)^2$, i.e. for $v\geq 0$,
\begin{eqnarray*}
I_S (v)^2 =4G^{-1}(v)G'(G^{-1}(v))^2
\end{eqnarray*}
is the restriction of a smooth function. 

Differentiating equation (\ref{VIS}) yields
\begin{eqnarray*}
V''=V'I'_S (V)=I_S(V)I'_S(V)=\frac{1}{2}\frac{dI_S^2}{dv}(V),
\end{eqnarray*}
The first derivative of $I_{S}$ is given by
\begin{eqnarray*}
\frac{d I_S}{dv}(V)=\frac{V''}{V'}=\frac{f'}{f},
\end{eqnarray*}
thus $I$ is convex if and only if $f'/f$ is nondecreasing.
Also $\displaystyle \frac{dI_S^2}{dv}(0)=2V''(0)=4\pi$.

Differentiating once more yields
\begin{eqnarray*}
V'''=\frac{1}{2}V'\frac{d^2 I_S^2}{dv^2}(V).
\end{eqnarray*}
The curvature, as a function of the distance $r$ to the origin, is given by
\begin{eqnarray*}
K(r)
&=&
-\frac{f''}{f}\\
&=&
-\frac{V'''}{V'}\\
&=&
-\frac{1}{2}\frac{d^2 I_S^2}{dv^2}(V(r))\\
&=&
-(I_{S}\frac{d^2 I_S}{dv^2}+(\frac{dI_S}{dv})^{2})(V(r)),
\end{eqnarray*}
which yields $\displaystyle \frac{d^2 I_S^2}{dv^2}(0)=-2K(0)$.

Conversely, let $\alpha$ be large. Choose first a smooth convex function $I_{0}$ on $\R_+$ such that $I_{0}(v)=v$ for $v\leq \alpha$, and $I'_{0}$ tends to $+\infty$. The corresponding surface of revolution satisfies $f(r)=e^{r}$ on $(-\infty,\log(\alpha/2\pi)]$, $f$ is nondecreasing, $f'/f$ is nondecreasing and tends to $+\infty$, and the curvature $-\frac{f''}{f}$ tends to $-\infty$. 

Given $\delta<\alpha$ and $k$ large enough, there exists a function $I=I_{\delta,k,\alpha}$ on $\R_+$ such that
\begin{itemize}
\item $I^{2}$ extends to a smooth function on a neighborhood of $\R_{+}$;
\item $I(0)=0$, $\displaystyle \frac{dI^2}{dv}(0)=4\pi$ and $\displaystyle \frac{d^2 I^2}{dv^2}(0)=-2k$;
\item $I(v)=I_{0}(v)$ for $v\geq\delta$;
\item $\displaystyle \frac{d^2 I^2}{dv^2}\geq -2k$ everywhere;
\item $I$ is nondecreasing and $v\mapsto I(v)/v$ is nondecreasing;
\item as $\delta$ tends to 0 (and $k$ tends to $+\infty$), $I=I_{\delta,k,\alpha}$ converges to $I_{0}$, $C^{\infty}$-uniformly on compact subsets of $\R_{-}$.
\end{itemize}
The inequality $I(v)\geq \sqrt{4\pi v -k\,v^2}$ for all $v$ follows from the second derivative bound.

\begin{figure}
\begin{center}
\includegraphics[width=2in]{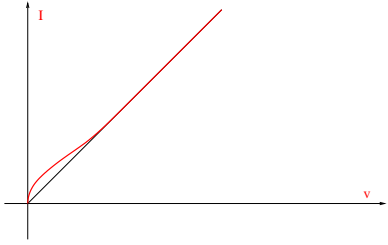}
\caption{Prescribed isoperimetric profile}
\label{I}
\end{center}
\end{figure}

The corresponding surface of revolution $S_{\delta,k,\alpha}$ is smooth, it satisfies $g_{\delta}=f$ on $[\log(\delta/2\pi),+\infty)$, $g_{\delta}$ is nondecreasing. Proposition \ref{profcusp} implies that for $\delta$ small enough, disks of revolution are extremal in $S_{\delta,k,\alpha}$, and the isoperimetric profile of $S_{\delta,k,\alpha}$ is equal to $I_{\delta,k,\alpha}$.\qed

\section{Cusp filling}
\label{filling}

Let $M$ be an orientable surface with at least two ends. Remove two points from $M$ and equip $M$ with a metric $g$ and calibration in such a way that neighborhoods of the deleted points be isometric to constant curvature $-1$ cusps. According to preceding sections, the isoperimetric profile of $(M,g)$ is linear. Furthermore, away from the cusp neighborhoods, $(M,g)$ has bounded curvature and injectivity radius. Let us cut cusp neighborhoods along horocycles of equal length and fill them with smooth caps of revolution where the isoperimetric profile is achieved by disks of revolution. We show that the isoperimetric profile of the obtained metric is smooth.

\subsection{The gluing procedure}

The following elementary properties will be needed.

\begin{lem}
\label{sousadditif}
Let $f$ and $g$ be nonnegative functions on $\R_+$ such that $v\mapsto f(v)/v$ and $v\mapsto g(v)/v$ are nonincreasing. Then $\min\{f,g\}$ has the same property.
Assume furthermore that $f(\delta)=g(\delta)$ for some $\delta>0$. Then the function $h$ such that $h=f$ on $[0,\delta]$ and $h=g$ on $(\delta,+\infty)$ also shares the same property. Any such function is subadditive, i.e.
\begin{eqnarray*}
f(v+v')\leq f(v)+f(v')
\end{eqnarray*}
for all $v$, $v'\in\R_+$.
\end{lem}

\begin{lem}
\label{incomplete}
Let $M$ denote a 2-dimensional complete Riemannian manifold with curvature bounded from above, $K\leq k$. Assume that there exists $\epsilon>0$ such that $M_{\epsilon}=\{x\in M \,;\,\mathrm{inj}(x)<\epsilon\}$ consists of a disjoint union of annuli.

Let $I^{d}$ denote the modified isoperimetric profile in which the competing domains are restricted to disjoint unions of domains of one of the following two types,
\begin{enumerate}
\item disks;
\item domains having at least one boundary component which is homotopic neither to zero nor to a component of $\partial M_{\epsilon}$.
\end{enumerate}
Then there exists $\eta>0$, $\eta<2\pi/k$, such that $I^{d}_{M}(v)\geq\sqrt{4\pi v -kv^{2}}$ for $v<\eta$.
\end{lem}

\proof
First, $I^{d}$ tends to zero, since small disks have small boundary length. A simple closed curve $c$ of length $<\epsilon$ which is not homotopic to zero is contained in $M_{\epsilon}$. Then $c$ is homotopic to a boundary component of $M_{\epsilon}$ with one of the two possible orientations. Therefore domains of the second type do not contribute to $I^{d}(v)$ for $v$ small. For disks, the Bol-Fiala isoperimetric inequality applies, \cite{Bol}, \cite{Fiala}. For unions of disks, use the subadditivity of $v\mapsto \sqrt{4\pi v -kv^{2}}$, which follows from Lemma \ref{sousadditif}.\qed

\begin{lem}
\label{cap}
Let $m$ be an integer and $0<\delta<\alpha/m$. Let $S^j_{\delta,\alpha}$, $j=1,\ldots,m$, denote $m$ smooth surfaces with boundary, of revolution, with area $\alpha$. Denote by $C^j_{\delta,\alpha}$ the rotationaly symmetric disk of revolution in $S^j_{\delta,\alpha}$ with area $\delta$. Assume that 
\begin{itemize}
  \item all $S^j_{\delta,\alpha}\setminus C^j_{\delta,\alpha}$ are constant curvature $-1$ annuli with inner boundary of length $\delta$;
  \item for all $v\leq m\delta$, $I_{S^j_{\delta,\alpha}}$ is achieved by a rotationally symmetric disk;
  \item $I_{S^j_{\delta,\alpha}}$ is nondecreasing, and $v\mapsto I_{S^j_{\delta,\alpha}}(v)/v$ is nonincreasing;
  \item for every $k$, there exists $\delta(k)>0$ such that for $v<\delta<\delta(k)$, $I_{S^j_{\delta,\alpha}}(v)\leq\sqrt{4\pi v -kv^{2}}$.
\end{itemize}

Let $M$ denote an orientable 2-dimensional complete Riemannian manifold. Assume that
\begin{itemize}
\item there exists $\epsilon\gg \delta$ such that $M_{\epsilon}=\{x\in M \,;\,\mathrm{inj}(x)<\epsilon\}$ consists of $m$ disjoint constant curvature $-1$ cusps;
\item let $\delta'$ be the number such that $\area(M_{\delta'})=m\delta$, then for every $v>m\delta$, $I_M (v)$ is achieved by a domain which contains $M_{\delta'}$;
\item $M$ has curvature bounded from above, $K\leq k$;
\item $I_M$ is nondecreasing;
\item for all $j=1,\ldots,m$, $I_M =I_{S^j_{\delta,\alpha}}$ on $[m\delta,m\alpha]$ where $\alpha\leq\area(M_{\epsilon})/m$ does not depend on $\alpha$.
\end{itemize}
Remove each component of $M_{\delta'}$ and replace it with a disk $C^j_{\delta,\alpha}$, to produce a smooth surface $N$. If $\delta$ is small enough, then $I_N$ is achieved, $I_{N}(v)=\min_j I_{S^j_{\delta,\alpha}}(v)$ for $v\leq m\delta$, $I_{N}(v)=I_M (v)$ for $v\geq m\delta$.
\end{lem}

\proof
Since the candidate isoperimetric profile is nondecreasing, when establishing the isoperimetric inequality for a domain $D$ in $N$, we can replace $D$ with any other domain with larger area and shorter boundary. For instance, if some boundary component $c$ of $D$ is null homotopic, it bounds a disk $D'$. Replace $D$ with $D\cup D'$. This increases the area of $D$ and decreases its boundary length. Therefore, we do not loose generality in assuming that $D$ is a disjoint union of domains which, unless they are disks, have no null homotopic boundary component. 

Denote by $N_{\epsilon}$ the complement of $M\setminus M_{\epsilon}$ in $N$. Let $w:N_{\epsilon}\to \R$ be the function, extending the 
(exponential of minus) horofunction on the constant curvature $-1$ annuli of $N$, such that $\area(\{w<t\})=t$. Then $\{w<m\delta\}$ is the union of the caps $C^j_{\delta,\alpha}$, $\{w<m\alpha\}$ is isometric to $\coprod_j S^j_{\delta,\alpha}$. If $\lng(\partial D\cap \{w<m\alpha\})\geq \lng(\{w=m\alpha\})$, then replacing $D$ with $D\cup\{w<m\alpha\}$ increases area and decreases boundary length. Thus we can assume that $\lng(\partial D\cap \{w<m\alpha\})< \lng(\{w=m\alpha\})$.

The distance in $M$ (and therefore in $N$) between the level sets $\{w=m^2 \delta\}$ and $\{w=m\alpha\}$ tends to $+\infty$ as $\delta$ tends to 0. If $\delta$ is small enough, there exists $t\in(m^2 \delta,m\alpha)$ such that $\partial D\cap\{w=t\}=\emptyset$. If the $j$-th component of $\{w=t\}$ is contained in $D$, then the whole $j$-th component of $\{w<t\}$ is contained in $D$. Otherwise, let us replace the part of $D$ contained in the $j$-th component of $\{w<t\}$ with the disk of the form $\{w<t_j \}$ in the same component which has the same area. Due to the isoperimetric inequality in $S^j_{\delta,\alpha}$, this decreases boundary length. Assume that $\area(D)<m\delta$. Then no component of $D$ can contain a component of $\{w<m^2 \delta\}$. Let $v_j$ denote the area of the part of $D$ contained in the $j$-th component of $\{w<t\}$ and $v'=\area(D)-\sum_j v_j =\area(D\cap\{w>t\})$. Note that Lemma \ref{incomplete} applies to $D\cap\{w>t\}$. We can assume that $m\delta<\eta$, the constant in this lemma, and that $\delta<\delta(k)$. Then
\begin{eqnarray*}
\lng(\partial D)
&\geq&
 \sum_j I_{S^j_{\delta,\alpha}}(v_j )+I^d_M(v')\\
&\geq&
\sum_j \min \{I_{S^j_{\delta,\alpha}}\}(v_j )+\min \{I_{S^j_{\delta,\alpha}}\}(v')\\
&\geq&
\min \{I_{S^j_{\delta,\alpha}}\}(\sum_j v_j +v')\\
&\geq&
\min \{I_{S^j_{\delta,\alpha}}\}(\area(D)),
\end{eqnarray*}
since $I^d_M (v)\geq \sqrt{4\pi v -k\,v^2}\geq\min \{I_{S^j_{\delta,\alpha}}\}$ for $v\in[0,m\delta]$ by Lemma \ref{incomplete} and $\min \{I_{S^j_{\delta,\alpha}}\}$ is subadditive, by Lemma \ref{sousadditif}.

If $\area(D)>m\delta$, let us construct a domain $D'$ in $M$ with the same area as $D$. The components of $D$ which do not intersect $\{w<t\}$ naturally live in $M$. The components of $D$ which contain components of $\{w<t\}$ also have obvious counterparts in $M$. To those components of $D$ which are contained in $\{w<t\}$, we associate a domain of the form $\{w<t_j \}$ in the corresponding end of $M$ with equal area. Since $I_M \leq I_{S^j_{\delta,\alpha}}$, the boundary length gets even smaller in $M$ than in $N$. This shows that $\lng(\partial D)\geq I_M (\area(D))$.

We have shown that $I_N \geq \min \{I_{S^j_{\delta,\alpha}}\}$ on $[0,m\delta]$ and $I_N \geq I_M$ elsewhere. Conversely, let $v\leq m\delta$. Let $j$ be the index for which $I_{S^j_{\delta,\alpha}}(v)$ is minimum. Then there exists in $S^j_{\delta,\alpha}$ a rotationally symmetric domain $D$ of area $v$ and minimal boundary length. This domains naturally embeds in $N$. Let $v>m\delta$. By assumption, $I_M (v)$ is achieved by a domain which contains $M_{\delta'}$. Again, this domain has an obvious counterpart in $N$ with the same area and the same boundary length. We conclude that $I_N$ is achieved and is equal to $\min \{I_{S^j_{\delta,\alpha}}\}$ on $[0,m\delta]$ and to $I_M$ elsewhere.\qed

\subsection{Smooth profile for 2 ended surfaces}
\label{proof2ends}

Let $\bar{M}$ be an orientable 2-manifold with at least two ends or one end of infinite genus. Remove two points from $\bar{M}$ to get $M$. Then $M$ satisfies the assumptions of Theorem \ref{main}. One can equip $M$ with a Morse function $w:M\to (0,+\infty)$ with compact level sets and a metric modelled on constant curvature $-1$ annuli, such that, for some $t>0$, $\{w<t\}$ consists of two cusps. Note that on $\{w\geq t\}$, the injectivity radius in bounded below and the curvature is bounded above. Theorem \ref{main} yields a calibration $\omega$ which garantees that the isoperimetric profile is linear $I_M (v)=v$.

Choose $\delta$ very small, according to Lemma \ref{cap}. Apply Corollary \ref{prescrcaps} in order to choose appropriate caps of revolution. Arrange so that one of the caps has an isoperimetric profile which is everywhere smaller than the other. Lemma \ref{cap} garantees that the resulting surface $N$ has a smooth isoperimetric profile. This completes the proof of the first part of Theorem \ref{2ends}.

\section{Construction of ultrahyperbolic surfaces }
\label{rearrangement}

\subsection{The rearrangement argument}

We need a generalization of Lemma \ref{calib}.

\begin{prop}
\label{rear}
Let $M$ be an oriented Riemannian manifold with volume element $vol_{M}$, let $\omega$ be an $n-1$-form on $M$ such that 
\begin{displaymath}
|\omega|\leq 1, \quad\textrm{ and }\quad d\omega=u\,vol_{M}
\end{displaymath}
for some nonnegative function $u$.
Let $V\subset M$ be a submanifold with compact boundary and finite volume $v$. Assume that
\begin{itemize}
  \item $\omega$ calibrates $\partial V$, i.e.
\begin{displaymath}
\omega_{|\partial V}=vol_{\partial V};
\end{displaymath}
  \item there exists $t$ such that $\{u<t\}\subset V\subset \{u\leq t\}$;
  \item there exist compact domains $V_{\epsilon}\subset V$ such that
\begin{displaymath}
vol(V\setminus V_{\epsilon})\to 0,\quad vol(\partial V \triangle \partial V_{\epsilon})\to 0
\end{displaymath}
as $\epsilon$ tends to 0.
\end{itemize}
Then $I_{M}(v)=\vol(\partial V)$.
\end{prop}

\proof
Let $D$ be a domain such that $\vol(D)=v=\vol(V)$. Then
\begin{eqnarray*}
\vol(V\setminus D)
&=&
\vol(V) -\vol(V\cap D)\\
&=&
\vol(D) -\vol(V\cap D)\\
&=&
\vol(D\setminus V).
\end{eqnarray*}

Since $u\leq t$ on $V$ and $u\geq t$ on $M\setminus V$,
\begin{eqnarray*}
\int_{V}u
&=&
\int_{V\cap D}u+\int_{V\setminus D}u\\
&\leq&
\int_{V\cap D}u+t\vol(V\setminus D)\\
&=&
\int_{V\cap D}u+t\vol(D\setminus V)\\
&\leq&
\int_{V\cap D}u+\int_{D\setminus V}u\\
&=&
\int_{D}u\\
&=&
\int_{D}d\omega\\
&=&
\int_{\partial D}\omega\\
&\leq&
\vol(\partial D),
\end{eqnarray*}
since $|\omega|\leq 1$. This shows that $I_{M}(v)\geq \int_{V}u$. 

Since $0\leq u\leq t$ on $V$, 
\begin{eqnarray*}
\int_{V\setminus V_{\epsilon}}u\leq t\vol(V\setminus V_{\epsilon})
\end{eqnarray*}
tends to 0. Also
\begin{eqnarray*}
|\int_{\partial V_{\epsilon}}\omega-\int_{\partial V}\omega|\leq\vol(\partial V \triangle \partial V_{\epsilon})
\end{eqnarray*}
tends to 0. Therefore
\begin{eqnarray*}
\int_{V}u=\lim_{\epsilon\to 0}\int_{V_{\epsilon}}u=\lim_{\epsilon\to 0}\int_{\partial V_{\epsilon}}\omega=\int_{\partial V}\omega.
\end{eqnarray*}
Since $\omega$ calibrates $\partial V$, $\int_{\partial V}\omega=vol(\partial V)$. We have shown that $I_{M}(v)\geq vol(\partial V)$, and that equality holds asymptotically for the domains $V_{\epsilon}$. We conclude that $I_{M}(v)=vol(\partial V)$.\qed

\subsection{Conformal changes of metrics}
\label{conf}

\begin{prop}
\label{conforme}
Let $(M,g)$ be an oriented Riemannian manifold, $w$ a smooth positive function on $M$, $\omega$ a calibration on $M$, i.e. $|\omega|=1$ and $d\omega=vol$. Assume that there exist $t_1 < t_0$ such that
\begin{itemize}
\item all critical values of $w$ belong to $(-\infty,t_{0})$;
\item $|\nabla w|=w$ on $\{w>t_{0}\}$;
\item on $\{w>t_{0}\}$, $\omega=\iota_{w^{-1}|\nabla w|}vol$;
\item for all $t\geq t_{0}$ or $t\leq t_1$, $\vol(\{w<t\})=t=\vol(\{w=t\})$;
\item for all $v< t_0$, $I_{M}(v)=v$ is asymptotically achieved by a sequence of domains which are contained in $\{w\leq t_0\}$.
\end{itemize}
Let $I$ be a smooth convex function on $\R_{+}$ such that $I(v)=v$ for $v\leq t_{0}$. Then there exists a conformal metric $g'=f^{2}g$ on some sublevel set $M'=\{w<\tau\}\subset M$ diffeomorphic to $M$ such that $g'\equiv g$ on $\{w\leq t_{0}\}$, $\vol(M')=+\infty$ and $I_{M'}=I$. 
\end{prop}

\proof
Let $\tau$ be some positive number, and $f:[0,\tau[\to \R_+$ a smooth function such that $f\equiv 1$ on $[0,t_0 ]$. Let $g'$ denote the conformal metric $f(w)^2 g$ on $M'=\{w<\tau\}\subset M$. Let $V(t)=\vol'(\{w<t\})$ denote the volumes of sublevel sets of $w$ in the new metric. According to the coarea formula, for $t\geq t_0$,
\begin{eqnarray*}
V(t)
&=&
V(t_0 )+\int_{\{t_0 <w<t\}}f(w)^n\\
&=&
t_0 +\int_{t_0}^{t}f(s)^n (\int_{\{w=s\}}\frac{1}{|\nabla w|})\,ds\\
&=&
t_0 +\int_{t_0}^{t}f(s)^n \,ds\\
&=&
\int_{0}^{t}f(s)^n \,ds,
\end{eqnarray*}
since $\vol(\{w=s\})=s=w=|\nabla w|$ on $\{w=s\}$ when $s\geq t_0$. In particular, $f(t)=V'(t)^{1/n}$.

Similarly, if $t\geq t_0$,
\begin{eqnarray*}
\vol'(\{w=t\})
&=&
f(t)^{n-1}\vol(\{w=t\})\\
&=&f(t)^{n-1}t.
\end{eqnarray*}
Assume that for all $v\geq t_0$, the isoperimetric profile $I=I_{(M',g')}$ is achieved by sublevel sets of $w$. Then, for all $t\geq t_0$,
\begin{eqnarray*}
\vol'(\{w=t\})=I(\vol'(\{w<t\})).
\end{eqnarray*}
This implies that
\begin{eqnarray*}
f(t)^{n-1}t=I(\int_{0}^{t}f(s)^n \,ds),
\end{eqnarray*}
i.e.
\begin{eqnarray}
V'(t)^{\frac{n-1}{n}}t=I(V(t)).\label{odeV}
\end{eqnarray}
This is a differential equation which, together with the initial condition $V(t_0 )=t_0$, uniquely determines $V$, and therefore $f=V'^{1/n}$.

Conversely, given a smooth function $I:\R_+ \to \R_+$ such that $I(v)=v$ for $v\leq t_0$, let $V:[0,\tau[\to \R_+$ be the maximal solution of Equation \ref{odeV} with initial condition $V(t_0 )=t_0$. Then $V(t)=t$ for $t\leq t_0$ and $\lim_{t\to\tau}V(t)=+\infty$. One sets $f=V'^{1/n}$. 

Let $\xi=\frac{\nabla w}{|\nabla w|}$. By assumption, the given calibration $\omega$ is equal to $\iota_{\xi}vol$ on $\{w\geq t_0 \}$. Let $\xi'=\frac{\nabla w}{|\nabla w|'}$ and $\omega'=\iota_{\xi'}vol'$. By definition, $|\omega'|'\leq 1$ and $\omega'$ calibrates level sets $\{w=t\}$ for $t\geq t_0$. Since
\begin{eqnarray*}
|\nabla w|'=f(w)^{-1}|\nabla w|,\quad \mathrm{and}\quad vol'=f(w)^n \,vol,
\end{eqnarray*}
\begin{eqnarray*}
\omega'=f(w)^{n-1}\omega,
\end{eqnarray*}
thus
\begin{eqnarray*}
d\omega'
&=&
(n-1)f(w)^{n-2}f'(w)\,dw\wedge\omega +f(w)^{n-1}d\omega\\
&=&
((n-1)f(w)^{n-2}f'(w)|\nabla w| +f(w)^{n-1})\,vol\\
&=&
u\,vol',
\end{eqnarray*}
where $u=h(w)$ and
\begin{eqnarray*}
 h(t)&=&
 f(t)^{-n}((n-1)f(t)^{n-2}f'(t)t+f(t)^{n-1})\\
 &=&
 (f^n )^{-1}\frac{d}{dt}(f^{n-1}t)\\
 &=&
 (\frac{dV}{dt})^{-1}\frac{dI(V(t))}{dt}\\
 &=&
 \frac{dI}{dv}(V(t)).
\end{eqnarray*}
If $I$ is convex, $h$ is nondecreasing, therefore large sublevel sets of $w$ satisfy
\begin{eqnarray*}
t\geq t_0 \quad\Rightarrow\quad \{u<h(t)\}\subset\{w<t\}\subset\{u\leq h(t)\}.
\end{eqnarray*}
Using, as compact approximations to $\{w< t\}$, the domains $\{\epsilon<w< t\}$, Lemma \ref{rear} applies, and for $v\geq t_0$, $I_{M'}(v)$ is achieved by the sublevel set of $w$ of volume $v$. For all $v$, $I_{M'}(v)\geq v$. By assumption, there exists for each $v<t_0$ a sequence of approximate extremal domains of area $v$ contained in $\{w\leq t_0 \}$. These domains are unaffected by the conformal change. Therefore $I_{M'}(v)=v$ for $v<t_0 $. This proves that $I_{M'}=I$.\qed

\subsection{Curvature and injectivity radius}

Let $\bar{M}$ be a compact orientable 2-dimensional Riemannian manifold. Remove finitely many (but at least two) points from $\bar{M}$ to get $M$. Choose two more points $a$ and $b$. Apply Theorem \ref{main} to $M\setminus \{a,b\}$. Get a Morse function $w:M\setminus \{a,b\}\to(0,+\infty)$ with compact level sets, a metric $g$ modelled on constant curvature $-1$ annuli except near critical points of $w$. Let $t_0$ (resp. $t_1$) be slightly larger (resp. smaller) than the largest (resp. smallest) critical value. Then $\{w<t_1 \}$ consists of two constant curvature $-1$ cusps and $\{w>t_0 \}$ of finitely many constant curvature $-1$ anticusps. Choose some smooth convex function $I$ on $\R_+$ such that $I(v)=v$ for $v\leq t_0$ and $I(v)=v\log(v)$ for large $v$. Apply Proposition \ref{conforme} to get a conformal metric on $M'=\{w<\tau\}$ with isoperimetric profile equal to $I$.

As observed in \cite{BP}, the second variation formula relates the first and second derivatives of the isoperimetric profile at $v$ to the Ricci curvature in the normal direction along the boundary of an extremal domain of volume $v$. Here, since the Gauss curvature $K$ is constant along the boundary of the extremal domain of volume $v$, one gets
\begin{eqnarray*}
\frac{d^{2}I_{M'}^{2}}{dv^{2}}(v)=-2K.
\end{eqnarray*}
Therefore, $K$ tends to $-\infty$ as $w$ tends to $\tau$.

Recall that the conformal change of metric is $g'=f(w)^2 g$ where $f=V'^{1/2}$ and $V$ is the solution of $V'(t)^{\frac{1}{2}}t=I(V(t))$ with initial condition $V(t_0 )=t_0$. The formulae
\begin{eqnarray*}
\frac{1}{t}-\frac{1}{\tau}=\int_{V(t)}^{+\infty}I(s)^{-2}ds
\end{eqnarray*}
and
\begin{eqnarray*}
f(t)=V'(t)^{1/2}=\frac{I(V(t))}{t}
\end{eqnarray*}
show that the asymptotic behaviour of $f(t)$ as $t$ tends to $\tau$ depends only on the asymptotic behaviour of $I(v)$ as $v$ tends to $+\infty$. An explicit calculation with $I(v)=v\log(v)$ shows that $\frac{1}{t}-\frac{1}{\tau}\sim V(t)^{-1}(\log(V(t)))^{-2}$, $\log(V(t))\sim\log(\tau t/\tau-t)$ and
\begin{eqnarray*}
f(t)\sim\frac{\tau t}{\tau-t}(\log(\frac{\tau t}{\tau-t}))^{-1}.
\end{eqnarray*}
In particular, $\int_{t_0}^{\tau}\frac{f(t)}{t}\,dt=+\infty$. This implies that the metric $g''=g+\frac{f(w)^2}{w^2}\,dw^2$ on $\{w<\tau\}$ is complete. Since
\begin{eqnarray*}
g'\geq g+\frac{f(w)^2 -1}{w^2}\,dw^2 ,
\end{eqnarray*}
the identity map $(\{f(w)\geq 2 \},g')\to (\{f(w)\geq 2 \},g'')$ is Lipschitz. This proves that $(M',g')$ is complete.

Let $x_j$ be a sequence of points in $M'$ such that $w(x_j)$ tends to $\tau$. Assume that $inj(x_j )$ is bounded by $L$. Since $M'$ is complete and negatively curved away from some sublevel set of $w$, for $j$ large, $inj(x_j )$ is equal to half the length of a geodesic loop $\ell_j$ which not null homotopic. Since $M$ has injectivity radius bounded below away from $\{w<t_1 \}$, the $g$-length of $\ell_j$ is bounded below. Since $\ell_j \subset B'(x_j ,L)\subset\{w>w(x_j )-\epsilon_j\}$ where $\epsilon_j$ tends to 0, $f$ is large at each point of $\ell_j$, a contradiction. We conclude that $inj(x)$ tends to $+\infty$ as $w(x)$ tends to $\tau$.

To complete the proof of Theorem \ref{2ends}, there remains to fill the cusps as described in paragraph \ref{proof2ends}.\qed

\begin{rem}
\label{maxcurv}
In the examples produced in Theorem \ref{2ends}, the isoperimetric profiles passes through three different regimes.
\begin{enumerate}
  \item At small areas, the profile is subadditive. One can arrange that the profile be achieved by a unique extremal domain. For this, it suffices to arrange that one of the caps be more positively curved than the other. There is some flexibility in the choice of the prescribed profile.
  \item At medium areas, the profile is exactly linear. There is no flexibility. The profile is achieved by large families of extremal domains.
  \item At large areas, the profile is superadditive. There is again flexibility in the choice of the profile. It is achieved by a unique extremal domain.
\end{enumerate}
\end{rem}

\section{Appendix : ultrahyperbolicity}

\begin{prop}
\label{achieved}
Let $M$ be a non compact complete 2-dimensional Riemannian manifold $M$ which is ultrahyperbolic, i.e.
\begin{itemize}
  \item curvature tends to $-\infty$ at infinity;
  \item injectivity radius tends to $+\infty$ at infinity.
\end{itemize}
On such a surface, the isoperimetric profile is achieved. Therefore, it is continuous on $\R_{+}$. Furthermore, for every $V<\area(M)$, the set of domains $D$ such that $\area(D)\in [0,V]$ and $\lng(\partial D)=I_{M}(\area(D))$, is compact.
\end{prop}

\proof
Given $\rho<0$ and $L>0$, let $C(\rho,L)$ denote the tubular neighborhood of width $2L$ of the set of points $x\in M$ such that either $inj(x)\leq L$ or $K(x)\geq \rho$. By assumption, this is a compact set. Denote by $J$ the function defined on $\R_{+}$ by $J_{\rho}(v)=\sqrt{4\pi v-\rho v^2}$.

We first show that for every $a>0$, if $\rho<-(\frac{4}{L}\cosh^{-1}(1+\frac{a}{2\pi}))^{2}$, then every connected domain $D$ with area $\leq a$ and boundary length $\leq L$ which is not contained in $C(\rho,L)$ satisfies the isoperimetric inequality $\lng(\partial D)\geq J_{\rho}(\area(D))$. Let $x\in D\setminus C(\rho,L)$. Then on the ball $B(x,2L)$, the injectivity radius is everywhere $\geq L$ and the curvature everywhere $\leq \rho$. Let $y\in B(x,L)$. Since $\area(B(y,L/4))\geq 2\pi(\cosh(\sqrt{-\rho} L/4)-1)>a$, the ball $B(y,L/4)$ overlaps the boundary of $D$. This implies that $d(x,y)< \frac{L}{4}+\frac{L}{4}+\frac{1}{2}\lng(\partial D)\leq L$. This show that $D$ is contained in $B(x,L)$. The Bol-Fiala inequality for disks (\cite{Bol}, \cite{Fiala}) implies that
\begin{eqnarray*}
\lng(\partial D)\geq J_{\rho}(\area(D)).
\end{eqnarray*}
Indeed, $D$ is homeomorphic to a disk with holes, filling the holes yields a disk $D'$, to which the Bol-Fiala isoperimetric inequality applies :
\begin{eqnarray*}
\lng(\partial D)
&\geq& 
\lng(\partial D')\\
&\geq&
J_{\rho}(\area(D'))\\
&\geq&
J_{\rho}(\area(D)),
\end{eqnarray*}
since $J_{\rho}$ is nondecreasing. Note that this implies that $\sqrt{-\rho}\area(D)\leq\lng(\partial D)$.

Let $D_{\ell}$ be a sequence of domains such that $\area(D_{\ell})$ converges to $v$, staying less than $a$, and $\lng(\partial D_{\ell})$ converges to $I_{M}(v)$, staying less than $L$. According to the compactness theorem for integral currents \cite{Federer}, one can diagonally extract a subsequence such that for every $\rho<0$, $D_{\ell}\cap C(\rho,L)$ converges in flat norm. Since $\area(\partial D_{\ell}\setminus C(\rho,L))\leq L/\sqrt{-\rho}$, the limiting current with unbounded support $D$ has area $v$. By semi-continuity, $\lng(\partial D)\leq I_{M}(v)$. Since $D$ minimizes boundary length for compactly supported area preserving perturbations, $D$ is a locally finite union of smooth domains, \cite{Fleming}.

Fix a point $x_0 \in M\setminus D_{\infty}$. Choose $\rho'<K(x_{0})$. When $r$ tends to 0, $\lng(\partial B(x_0 ,r))$ and $\area(B(x_{0},r))$ behave asymptotically like in constant curvature $K(x_{0})$. It follows that, for $w$ small enough (say $w\leq v_0$), the function $I$ such that $\lng(\partial B(x_0 ,r))=I(\area(B(x_{0},r))$ satisfies $I(w)\leq\sqrt{4\pi w -\rho'w^2}$. Let $v_{1}=\min\{v_{0},\area(B(x_{0},dist(x_{0},D)))\}$. Given $a$ and $L$, choose $\rho<\rho'$ such that $\rho<-(\frac{4}{L}\cosh^{-1}(1+\frac{a}{2\pi}))^{2}$ and $\frac{L}{\sqrt{-\rho}}<v_1$. Let $B=B(x_0 ,r_0 )$ denote the ball centered at $x_0$ with area $v_1$. By construction, $B$ does not intersect $D$.

Assume that $D$ is not contained in $C(\rho,L)$. Let $D_{2}$ be the union of components of $D$ which are not contained in $C(\rho,L)$, and $D_{1}=D\setminus D_{2}$. For $r>0$, let $D_{r}=D_{1}\cup B(x_{0},r)$. Since every component $\Delta$ of $D$ which is not contained in $C(\rho,L)$ satisfies $\sqrt{-\rho}\area(\Delta)\leq\lng(\partial \Delta)$, $\area(D_2 )\leq \frac{L}{\sqrt{-\rho}}<v_0$. Therefore there exists $r\in(0,r_0)$ such that $\area(D_r )=\area(D)$. Then
\begin{eqnarray*}
\lng(\partial D_r )-\lng(\partial D_1 )
&=&
\lng(\partial B(x_0 ,r))\\
&\leq&
I(\area(B(x_0 ,r)))\\
&\leq&
\sqrt{4\pi \alpha-k\alpha^2},
\end{eqnarray*}
where $\alpha=\area(D_2 )$. Since $v\mapsto J_{\rho}(v)/v$ is nonincreasing, $J_{\rho}$ is subadditive. Therefore the inequality $\lng(\partial \Delta)\geq J_{\rho}(\area(\Delta))$, valid for each component $\Delta$ of $D_2$, implies $\lng(\partial D_2 )\geq J_{\rho}(\area(D_2 ))$. It follows that $\lng(\partial D_r )<\lng(\partial D)$. This contradicts the extremality of $D$. We conclude that $D$ is compact.

This proves the existence of compact extremal domains. According to \cite{BP}, this implies that the isoperimetric profile is continuous. The compactness of the set of extremal domains follows.\qed

\begin{rem}
\label{endultra}
Here is a relative version of Proposition \ref{achieved}. \end{rem}
Let $M$ be a 2-dimensional Riemannian manifold. Say an end $E$ of $M$ is \emph{ultrahyperbolic} if the injectivity radius tends to $+\infty$ and curvature to $-\infty$ in $E$. Then for every minimizing sequence $D_{\ell}$ for the isoperimetric problem in $M$, there is a subsequence which does not enter in $E$.

Relative ultrahyperbolicity is used in Proposition \ref{profcusp}.

\bibliographystyle{plain}
\bibliography{calibration}

\vskip1cm
\noindent
Renata Grimaldi\\
Universit\`a di Palermo\\
Dipartimento di Metodi e Modelli Matematici\\
Facolta di Ingegneria\\
Vialle delle Scienze - 90128 Palermo (Italia)\\
\smallskip\noindent
{\tt grimaldi@unipa.it}
\par\medskip\noindent
Pierre Pansu\\
Laboratoire de Math\'ematique d'Orsay\\
UMR 8628 du C.N.R.S.\\
B\^atiment 425\\
Universit\'e Paris-Sud 11 - 91405 Orsay (France)\\
\smallskip\noindent
{\tt\small Pierre.Pansu@math.u-psud.fr}\\
http://www.math.u-psud.fr/$\sim$pansu

\end{document}